\documentclass[final]{siamltex}
\usepackage{amsfonts,epsfig,multicol}

\newcommand{\diam}{\mathop{\mathrm{diam}}}
\newcommand{\bigO}{\mathcal{O}}
\newcommand{\trans}{\mathsf{T}}
\newcommand{\FMM}{\mathrm{FMM}}
\newcommand{\RS}{\mathrm{RS}}
\newcommand{\row}{\mathrm{r}}
\newcommand{\col}{\mathrm{c}}
\newcommand{\cm}{\mathrm{cm}}
\newcommand{\mv}{\mathrm{mv}}

\newcommand{\lu}{\mathrm{lu}}
\newcommand{\sv}{\mathrm{sv}}
\newcommand{\E}{{\tiny E}}

\title{A fast direct solver for structured linear systems by recursive skeletonization
\thanks{This work was supported in part by the National Science Foundation under grant
DGE-0333389, by the U.S.\ Department of Energy under contract DEFG0288ER25053, and by
the Air Force Office of Scientific Research under NSSEFF Program Award FA9550-10-1-0180.}}
\author{Kenneth L. Ho
\thanks{Courant Institute of Mathematical Sciences and Program in Computational Biology, New York University, New York, NY, USA ({\tt ho@courant.nyu.edu}).}
\and Leslie Greengard
\thanks{Courant Institute of Mathematical Sciences, New York University, New York, NY, USA ({\tt greengard@courant.nyu.edu}).}}

\begin{document}

 \maketitle

 \begin{abstract}
  We present a fast direct solver for structured linear systems based on multilevel matrix compression. Using the recently developed interpolative decomposition of a low-rank matrix in a recursive manner, we embed an approximation of the original matrix into a larger, but highly structured sparse one that allows fast factorization and application of the inverse. The algorithm extends the Martinsson/Rokhlin method developed for 2D boundary integral equations and proceeds in two phases: a precomputation phase, consisting of matrix compression and factorization, followed by a solution phase to apply the matrix inverse. For boundary integral equations which are not too oscillatory, e.g., based on the Green's functions for the Laplace or low-frequency Helmholtz equations, both phases typically have complexity $\bigO (N)$ in two dimensions, where $N$ is the number of discretization points. In our current implementation, the corresponding costs in three dimensions are $\bigO (N^{3/2})$ and $\bigO (N \log N)$ for precomputation and solution, respectively. Extensive numerical experiments show a speedup of $\sim 100$ for the solution phase over modern fast multipole methods; however, the cost of precomputation remains high. Thus, the solver is particularly suited to problems where large numbers of iterations would be required. Such is the case with ill-conditioned linear systems or when the same system is to be solved with multiple right-hand sides. Our algorithm is implemented in Fortran and freely available.
 \end{abstract}



 \pagestyle{myheadings}
 \thispagestyle{plain}
 \markboth{K.\ L.\ HO AND L.\ GREENGARD}{FAST DIRECT SOLVER BY RECURSIVE SKELETONIZATION}

 \section{Introduction}
 \label{sec:introduction}
 Many problems in computational science and engineering require the solution of large, dense linear systems. Standard direct methods based on Gaussian elimination, of course, require $\bigO (N^{3})$ work, where $N$ is the system size. This quickly becomes infeasible as $N$ increases. As a result, such systems are typically solved iteratively, combining GMRES \cite{saad:1986:siam-j-sci-stat-comput}, Bi-CGSTAB \cite{van-der-vorst:1992:siam-j-sci-stat-comput} or some other iterative scheme with fast algorithms to apply the system matrix, when available. For the integral equations of classical physics, this combination has led to some of the fastest solvers known today, with dramatically lower complexity estimates of the order $\bigO (N)$ or $\bigO (N \log N)$ \cite[and references therein]{chew:2001:artech-house,liu:2009:cambridge,nishimura:2002:appl-mech-rev}.

 Despite their tremendous success, however, iterative methods still have several significant disadvantages when compared with their direct counterparts:
 \begin{romannum}
  \item
   {\em The number of iterations required by an iterative solver is highly sensitive to the conditioning of the system matrix.} Ill-conditioning arises, for example, in the solution of problems near resonance (particularly in the high frequency regime), in geometries with ``close-to-touching'' interactions, in multi-component physics models with large contrasts in material properties, etc. Under these circumstances, the solution time can be far greater than expected. Direct methods, by contrast, are robust in the sense that their solution time does not degrade with conditioning. Thus, they are often preferred in production environments, where reliability of the solver and predictability of the solution time are important.
  \item
   {\em One often wishes to solve a linear system governed by a fixed matrix with multiple right-hand sides.} This occurs, for example, in scattering problems, in optimization, and in the modeling of time-dependent processes in fixed geometry. Most iterative methods are unable to effectively exploit the fact that the system matrix is the same, and simply treat each right-hand side as a new problem. Direct methods, on the other hand, are extremely efficient in this regard: once the system matrix has been factored, the matrix inverse can be applied to each right-hand side at a much lower cost.
  \item
   {\em One often wishes to solve problems when the system matrix is altered by a low-rank modification.} Standard iterative methods do a poor job of exploiting this fact. Direct methods, on the other hand, can update the factorization of the original matrix using the Sherman-Morrison-Woodbury formula \cite{hager:1989:siam-rev} or use the existing factorization as a preconditioner.
 \end{romannum}

 In this paper, we present an algorithm for the solution of structured linear systems that overcomes these deficiencies, while remaining competitive with modern fast iterative solvers in many practical situations. The algorithm directly constructs a compressed (``data-sparse'') representation of the system matrix inverse, assuming only that the matrix has a block low-rank structure similar to that utilized by fast matrix-vector product techniques like the fast multipole method (FMM) \cite{greengard:1987:j-comput-phys,greengard:1997:acta-numer}. Such matrices typically arise from the discretization of integral equations, where the low-rank structure can be understood in terms of far-field interactions between clusters of points, but the procedure is general and makes no {\it a priori} assumptions about rank. Our scheme is a multilevel extension of the work described in \cite{greengard:2009:acta-numer}, which itself is based on the fast direct multilevel method developed for 2D boundary integral equations by Martinsson and Rokhlin \cite{martinsson:2005:j-comput-phys}.

 While we do not seek to review the literature on fast direct solvers here, it is worth noting that similar efforts have been (and continue to be) pursued by various groups, most notably in the context of hierarchically semiseparable (HSS) matrices \cite{chandrasekaran:2006b:siam-j-matrix-anal-appl,chandrasekaran:2006a:siam-j-matrix-anal-appl,xia:2009:siam-j-matrix-anal-appl} and $\mathcal{H}$-matrices \cite{hackbusch:1999:computing,hackbusch:2002:computing,hackbusch:2000:computing}. A short historical discussion can be found in \cite{greengard:2009:acta-numer} as well as in the recent article by Gillman {\it et al.} \cite{gillman:2012:front-math-china}. The latter paper makes several improvements on the algorithm of \cite{martinsson:2005:j-comput-phys}, and presents a simple framework for understanding, implementing, and analyzing schemes for inverting integral equations on curves (that is, domains parametrized by a single variable). Planar domains with corners were treated recently in \cite{bremer:preprint}. Applications to electromagnetic wave problems were considered in \cite{wei:2011:27th-annual-rev-prog-appl-comput-electromag,winebrand:2009:intl-conf-electromag-adv-application}. Finally, it should be noted that Gillman's dissertation \cite{gillman:2011:univ-colorado-boulder} includes 3D experiments that also extend the Martinsson-Rokhlin formalism to the case of integral equations on surfaces.

 The present paper provides a mix of analysis, algorithmic work, and applications. The novelty of our contribution lies:
 \begin{romannum}
  \item
   in the use of compression and auxilliary variables to embed an approximation of the original dense matrix into a sparse matrix framework that can make use of standard and well-developed sparse matrix technology;
  \item
   in providing detailed numerical experiments in both 2D and 3D; and
  \item
   in demonstrating the utility of fast direct solvers in several applications.
 \end{romannum}
 We believe that the scheme is substantially simpler to implement than prior schemes and that it leads to a more stable solution process.

 As in previous schemes (see, e.g., \cite{gillman:2012:front-math-china}), the core algorithm in our work computes a compressed matrix representation using the interpolative decomposition (ID) \cite{cheng:2005:siam-j-sci-comput,liberty:2007:proc-natl-acad-sci-usa,woolfe:2008:appl-comput-harmon-anal} via a multilevel procedure that we refer to as {\em recursive skeletonization}. Once obtained, the compressed representation serves as a platform for fast matrix algebra including matrix-vector multiplication and matrix inversion. In its former capacity, the algorithm may be viewed as a generalized or kernel-independent FMM \cite{gimbutas:2002:siam-j-sci-comput,martinsson:2007:siam-j-sci-comput,ying:2004:j-comput-phys}; we explore this application in \S \ref{sec:numerical-examples}. For matrix inversion, we show how to embed the compressed representation in an equivalent (but larger) sparse system, much in the style of \cite{chandrasekaran:2006b:siam-j-matrix-anal-appl,pals:2004:univ-calif-santa-barbara}. We then use a state-of-the-art sparse matrix solver to do the rest. We are grateful to David Bindel for initially suggesting an investigation of the sparse matrix formalism and rely in this paper on the sparse direct solver software UMFPACK \cite{davis:2004:acm-trans-math-softw,davis:1997:siam-j-matrix-anal-appl}. As in dense LU factorization, the direct solver is a two-phase process. First, following the generation of the compressed matrix embedding, a factored representation of the inverse is constructed. Second, in the solution phase, the matrix inverse is applied in a rapid manner to a specified right-hand side. As expected, the solution phase is very inexpensive, often beating a single FMM call by several orders of magnitude. For boundary integral equations without highly oscillatory kernels, e.g., the Green's function for the Laplace or low-frequency Helmholtz equation, both phases typically have complexity $\bigO (N)$ in 2D. In 3D, the complexities in our current implementation are $\bigO (N^{3/2})$ and $\bigO (N \log N)$ for precomputation (compression and factorization) and solution, respectively.

 The remainder of this paper is organized as follows. In \S \ref{sec:preliminaries}, we define the matrix structure of interest and review certain aspects of the ID. In \S \ref{sec:algorithm}, we review the recursive skeletonization algorithm for matrix compression and describe the new formalism for embedding the compressed matrix in a sparse format. In \S \ref{sec:complexity-analysis}, we study the complexity of the algorithm for non-oscillatory problems, while in \S \ref{sec:error-analysis}, we give error estimates for applying a compressed matrix and its inverse. In \S \ref{sec:numerical-examples}, we demonstrate the efficiency and generality of our scheme by reporting numerical results from its use as a generalized FMM, as a direct solver, and as an accelerator for molecular electrostatics and scattering problems. Finally, in \S \ref{sec:generalizations-conclusions}, we summarize our findings and discuss future work.

 \section{Preliminaries}
 \label{sec:preliminaries}
 In this section, we discuss the precise matrix structure that makes our fast solver possible. For this, let $A \in \mathbb{C}^{N \times N}$ be a matrix whose index vector $J = (1, 2, \dots, N)$ is grouped into $p$ contiguous blocks of $n_{i}$ elements each, where $\sum_{i = 1}^{p} n_{i} = N$:
 $$
  J_{i} = \left( \sum_{j = 1}^{i - 1} n_{j} + 1, \sum_{j = 1}^{i - 1} n_{j} + 2, \dots, \sum_{j = 1}^{i} n_{j} \right), \quad i = 1, \dots, p.
 $$
 Then the linear system $A \mathbf{x} = \mathbf{b}$ can be written in the form
 $$
  \sum_{j = 1}^{p} A_{ij} \mathbf{x}_{j} = \mathbf{b}_{i}, \quad i = 1, \dots, p,
 $$
 where $\mathbf{x}_{i}, \mathbf{b}_{i} \in \mathbb{C}^{n_{i}}$ and $A_{ij} \in \mathbb{C}^{n_{i} \times n_{j}}$. Solution of the full linear system by classical Gaussian elimination is well-known to require $\bigO (N^{3})$ work.

 \begin{definition}[block separability]
  The matrix A is said to be {\em block separable} if each off-diagonal submatrix $A_{ij}$ can be decomposed as the product of three low-rank matrices:
  \begin{equation}
   A_{ij} = L_{i} S_{ij} R_{j}, \quad i \neq j,
   \label{eq:block separable}
  \end{equation}
  where $L_{i} \in \mathbb{C}^{n_{i} \times k^{\row}_{i}}$, $S_{ij} \in \mathbb{C}^{k^{\row}_{i} \times k^{\col}_{j}}$, and $R_{j} \in \mathbb{C}^{k^{\col}_{j} \times n_{j}}$, with $k^{\row}_{i}, k^{\col}_{i} \ll n_{i}$. Note that in (\ref{eq:block separable}), the left matrix $L_{i}$ depends only on the index $i$ and the right matrix $R_{j}$ depends only on the index $j$.
 \end{definition}

 We will see how such a factorization arises below. The term {\em block separable} was introduced in \cite{gillman:2012:front-math-china}, and is closely related to that of semiseparable matrices \cite{chandrasekaran:2006b:siam-j-matrix-anal-appl,chandrasekaran:2006a:siam-j-matrix-anal-appl,xia:2009:siam-j-matrix-anal-appl} and $\mathcal{H}$-matrices \cite{hackbusch:1999:computing,hackbusch:2002:computing,hackbusch:2000:computing}. In \cite{greengard:2009:acta-numer}, the term {\em structured} was used, but block separable is somewhat more informative.

 \begin{definition}[off-diagonal block rows and columns]
  The {\em $i$th off-diagonal block row} of $A$ is the submatrix $[ \, A_{i1} \ \cdots \ A_{i(i - 1)} \ A_{i(i + 1)} \ \cdots \ A_{ip} \, ]$ consisting of the $i$th block row of $A$ with the diagonal block $A_{ii}$ deleted; the {\em off-diagonal block columns} of $A$ are defined analogously.
 \end{definition}

 Clearly, the block separability condition (\ref{eq:block separable}) is equivalent to requiring that the $i$th off-diagonal block row and column have rank $k^{\row}_{i}$ and $k^{\col}_{i}$, respectively, for $i = 1, \dots, p$ (see \S \ref{sec:algorithm} for details).

 When $A$ is block separable, it can be written as
 \begin{equation}
  A = D + LSR,
  \label{eq:compressed-representation}
 \end{equation}
 where
 $$
  D = \left[
  \begin{array}{ccc}
   A_{11}\\
   & \ddots\\
   & & A_{pp}
  \end{array} \right] \in \mathbb{C}^{N \times N}
 $$
 is block diagonal, consisting of the diagonal blocks of $A$,
 $$
  L = \left[
  \begin{array}{ccc}
   L_{1}\\
   & \ddots\\
   & & L_{p}
  \end{array} \right] \in \mathbb{C}^{N \times K_{\row}}, \qquad R = \left[
  \begin{array}{ccc}
   R_{1}\\
   & \ddots\\
   & & R_{p}
  \end{array} \right] \in \mathbb{C}^{K_{\col} \times N}
 $$
 are block diagonal, where $K_{\row} = \sum_{i = 1}^{p} k^{\row}_{i}$ and $K_{\col} = \sum_{i = 1}^{p} k^{\col}_{i}$, and
 $$
  S = \left[
  \begin{array}{cccc}
   0 & S_{12} & \cdots & S_{1p}\\
   S_{21} & 0 & \cdots & S_{2p}\\
   \vdots & \vdots & \ddots & \vdots\\
   S_{p1} & S_{p2} & \cdots & 0\\
  \end{array} \right] \in \mathbb{C}^{K_{\row} \times K_{\col}}
 $$
 is dense with zero diagonal blocks. It is convenient to let $\mathbf{z} = R \mathbf{x}$ and $\mathbf{y} = S \mathbf{z}$. We can then write the original system in the form
 \begin{equation}
  \left[
  \begin{array}{ccc}
   D & L\\
   R & & -I\\
   & -I & S
  \end{array} \right] \left[
  \begin{array}{c}
   \mathbf{x}\\
   \mathbf{y}\\
   \mathbf{z}
  \end{array} \right] = \left[
  \begin{array}{c}
   \mathbf{b}\\
   \mathbf{0}\\
   \mathbf{0}
  \end{array} \right].
  \label{eq:sparse-embedding}
 \end{equation}
 This system is highly structured and sparse, and can be efficiently factored using standard techniques. If we assume that each block corresponds to $N_{i} = N/p$ unknowns and that the ranks $k^{\row}_{i} = k^{\col}_{i} \equiv k$ of the off-diagonal blocks are all the same, it is straightforward to see \cite{gillman:2012:front-math-china,greengard:2009:acta-numer} that a scheme based on (\ref{eq:compressed-representation}) or (\ref{eq:sparse-embedding}) requires an amount of work of the order $\bigO(p (N/p)^{3} + p^{3} k^{3})$.

 In many contexts (including integral equations), the notion of block separability is applicable on a hierarchy of subdivisions of the index vector. That is to say, a decomposition of the form (\ref{eq:compressed-representation}) can be constructed at each level of the hierarchy. When a matrix has this structure, much more powerful solvers can be developed, but they will require some additional ideas (and notation).

 \subsection{Hierarchically structured matrices}
 \label{sec:preliminaries:matrix-structure}
\begin{figure}
  \begin{center}
   \includegraphics{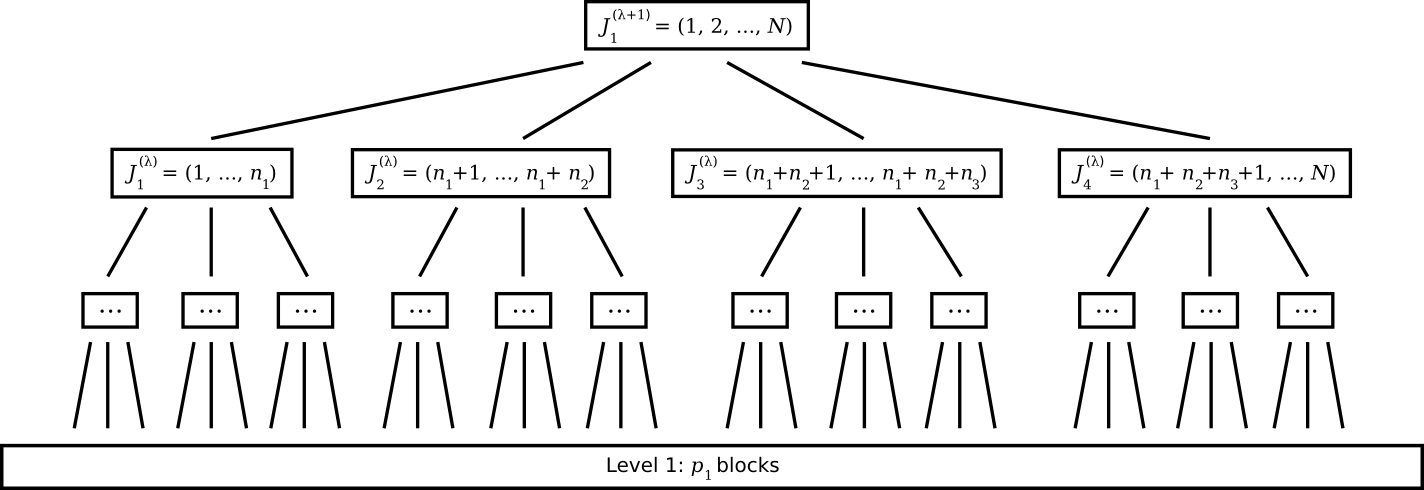}
  \end{center}
  \caption{An example of a tree structure imposed on the index vector $(1, 2, \dots, N)$. At each level of the hierarchy, a contiguous block of indices is divided into a set of children, each of which corresponds to a contiguous subset of indices.}
  \label{fig:tree-struct}
 \end{figure}

 Our treatment in this section follows that of \cite{gillman:2012:front-math-china}. Let $J = (1, 2, \dots, N)$ be the index vector of a matrix $A \in \mathbb{C}^{N \times N}$. We assume that a tree structure $\tau$ is imposed on $J$ which is $\lambda + 1$ levels deep. At level $l$, we assume that there are $p_{l}$ nodes, with each such node $J^{(l)}_{i}$ corresponding to a contiguous subsequence of $J$ such that
 $$
  \left\{ J^{\left( l \right)}_{1}, J^{\left( l \right)}_{2}, \dots, J^{\left( l \right)}_{p_{l}} \right\} = J.
 $$
 We denote the {\em finest level} as level $1$ and the coarsest level as level $\lambda + 1$ (which consists of a single block). Each node $J^{(l)}_{i}$ at level $l > 1$ has a finite number of children at level $l - 1$ whose concatenation yields the indices in $J^{(l)}_{i}$ (Fig.\ \ref{fig:tree-struct}).

\begin{figure}
  \begin{center}
   \includegraphics{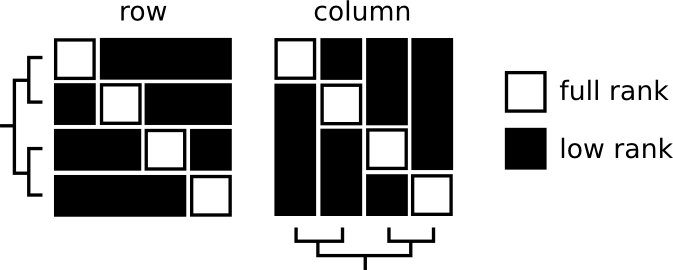}
  \end{center}
  \caption{Matrix rank structure. At each level of the index tree, the off-diagonal block rows and columns (black) must have low numerical rank; the diagonal blocks (white) can in general be full-rank.}
  \label{fig:mat-struct}
 \end{figure}

 The matrix $A$ is {\em hierarchically block separable} \cite{gillman:2012:front-math-china} if it is block separable at each level of the hierarchy defined by $\tau$. In other words, it is structured in the sense of the present paper if, on each level of $\tau$, the off-diagonal block rows and columns are low-rank (Fig.\ \ref{fig:mat-struct}).
 Such matrices arise, for example, when discretizing integral equations with non-oscillatory kernels (up to a specified precision).

 {\it Example 1}. Consider the integral operator
 \begin{equation}
  \phi \left( x \right) = \int G \left( x, y \right) \rho \left( y \right) dy
  \label{eq:integral-operator}
 \end{equation}
 where
 \begin{equation}
  G \left( x, y \right) = -\frac{1}{2 \pi} \log \left| x - y \right|
  \label{eq:2d-laplace}
\end{equation}
 is the Green's function for the 2D Laplace equation, and the domain of integration is a square $B$ in the plane. This is a 2D {\em volume integral operator}. Suppose now that we discretize (\ref{eq:integral-operator}) on a $\sqrt{N} \times \sqrt{N}$ grid:
 \begin{equation}
  \phi \left( x_{i} \right) = \frac{1}{N} \sum_{j \neq i}  G \left( x_{i}, x_{j} \right) \rho \left( x_{j} \right).
  \label{eq:integral-operator-discret}
 \end{equation}
 (This is not a high-order quadrature but that is really a separate issue.) Let us superimpose on $B$ a quadtree of depth $\lambda + 1$, where $B$ is the root node (level $\lambda + 1$). Level $\lambda$ is obtained from level $\lambda + 1$ by subdividing the box $B$ into four equal squares and reordering the points $x_{i}$ so that each child holds a contiguous set of indices. This procedure is carried out until level $1$ is reached, reordering the nodes at each level so that the points contained in every node at every level correspond to a contiguous set of indices. It is clear that, with this ordering, the matrix corresponding to (\ref{eq:integral-operator-discret}) is hierarchically block separable, since the interactions between nonadjacent boxes at every level are low-rank to any specified precision (from standard multipole estimates \cite{greengard:1987:j-comput-phys}). Adjacent boxes are low-rank for a more subtle reason (see \S \ref{sec:complexity-analysis} and Fig.\ \ref{fig:recur-subdiv}).

 {\it Example 2}. Suppose now that we wish to solve an interior Dirichlet problem for the Laplace equation in a simply connected 3D domain $\Omega$ with boundary $\partial \Omega$:
 \begin{equation}
  \Delta u = 0 \quad \mbox{in $\Omega$}, \qquad u = f \quad \mbox{on $\partial \Omega$}.
  \label{eq:laplace-dirichlet}
 \end{equation}
 Potential theory \cite{guenther:1988:prentice-hall} suggests that we seek a solution in the form of a double-layer potential
 \begin{equation}
  u \left( x \right) = \int_{\partial \Omega} \frac{\partial G}{\partial \nu_{y}} \left( x, y \right) \sigma \left( y \right) dy \quad \mbox{for $x \in \Omega$},
  \label{eq:double-layer}
 \end{equation}
 where
 \begin{equation}
  G \left( x, y \right) = \frac{1}{4 \pi \left| x - y \right|}
  \label{eq:3d-laplace}
 \end{equation}
 is the Green's function for the 3D Laplace equation, $\nu_{y}$ is the unit outer normal at $y \in \partial \Omega$, and $\sigma$ is an unknown surface density. Letting $x$ approach the boundary, this  gives rise to the second-kind Fredholm equation
 \begin{equation}
  -\frac{1}{2} \sigma \left( x \right) + \int_{\partial \Omega} \frac{\partial G}{\partial \nu_{y}} \left( x, y \right) \sigma \left( y \right) dy = f \left( x \right).
  \label{eq:2k-integral-eqn}
 \end{equation}
 Using a standard collocation scheme based on piecewise constant densities over a triangulated surface, we enclose $\partial \Omega$ in a box $B$ and bin sort the triangle centroids using an octree where, as in the previous example, we reorder the nodes so that each box in the hierarchy contains contiguous indices. It can be shown that the resulting matrix is also hierarchically block separable (see \S \ref{sec:complexity-analysis} and \cite{greengard:2009:acta-numer}).

 We turn now to a discussion of the ID, the compression algorithm that we will use to compute low-rank approximations of off-diagonal blocks. A useful feature of the ID is that it is able to compute the rank of a matrix on the fly, since the exact ranks of the blocks are difficult to ascertain {\em a priori}---that is to say, the ID is {\em rank-revealing}.

 \subsection{Interpolative decomposition}
 \label{sec:preliminaries:interpolative-decomposition}
 Many decompositions exist for low-rank matrix approximation, including the singular value decomposition, which is well-known to be optimal \cite{golub:1996:johns-hopkins-univ}. Here, we consider instead the ID \cite{cheng:2005:siam-j-sci-comput,liberty:2007:proc-natl-acad-sci-usa,woolfe:2008:appl-comput-harmon-anal}, which produces a near-optimal representation that is more useful for our purposes as it permits an efficient scheme for multilevel compression when used in a hierarchical setting.

 \begin{definition}[interpolative decomposition]
  Let $A \in \mathbb{C}^{m\times n}$ be a matrix, and $\| \cdot \|$ the matrix $2$-norm. A rank-$k$ approximation of $A$ in the form of an {\em interpolative decomposition (ID)} is a representation $A \approx BP$, where $B \in \mathbb{C}^{m \times k}$, whose columns constitute a subset of the columns of $A$, and $P \in \mathbb{C}^{k \times n}$, a subset of whose columns makes up the $k \times k$ identity matrix, such that $\left\| P \right\|$ is small and $\| A - BP \| \sim \sigma_{k + 1}$, where $\sigma_{k + 1}$ is the ($k + 1$)st greatest singular value of $A$. We call $B$ and $P$ the {\em skeleton} and {\em projection matrices}, respectively.
 \end{definition}

 Clearly, the ID compresses the column space of $A$; to compress the row space, simply apply the ID to $A^{\trans}$, which produces an analogous representation $A = \tilde{P} \tilde{B}$, where $\tilde{P} \in \mathbb{C}^{m \times k}$ and $\tilde{B} \in \mathbb{C}^{k \times n}$.

 \begin{definition}[row and column skeletons]
  The row indices that corrrespond to the retained rows in the ID are called the {\em row} or {\em incoming skeletons}. The column indices that corrrespond to the retained columns in the ID are called the {\em column} or {\em outgoing skeletons}.
 \end{definition}

 Reasonably efficient schemes for constructing an ID exist \cite{cheng:2005:siam-j-sci-comput,liberty:2007:proc-natl-acad-sci-usa,woolfe:2008:appl-comput-harmon-anal}. By combining such schemes with methods for estimating the approximation error, we can compute an ID to any relative precision $\epsilon > 0$ by adaptively determining the required rank $k$ \cite{liberty:2007:proc-natl-acad-sci-usa}. This is the sense in which we will use the ID.

 While previous related work \cite{greengard:2009:acta-numer,martinsson:2005:j-comput-phys} used the deterministic $\bigO (kmn)$ algorithm of \cite{cheng:2005:siam-j-sci-comput}, we employ here the latest compression technology based on random sampling, which typically requires only $\bigO (mn \log k + k^{2} n)$ operations \cite{liberty:2007:proc-natl-acad-sci-usa,woolfe:2008:appl-comput-harmon-anal}.

 \section{Algorithm}
 \label{sec:algorithm}
 In this section, we first describe the ``standard'' ID-based fast multilevel matrix compression algorithm (as in \cite{gillman:2012:front-math-china,martinsson:2005:j-comput-phys}). The HSS and $\mathcal{H}$-matrix formalisms use the same underlying philosophy \cite{chandrasekaran:2006b:siam-j-matrix-anal-appl,chandrasekaran:2006a:siam-j-matrix-anal-appl,hackbusch:1999:computing,hackbusch:2002:computing,hackbusch:2000:computing,xia:2009:siam-j-matrix-anal-appl}. We then describe our new inversion scheme.

 \subsection{Hierarchical matrix compression}
 Let $A \in \mathbb{C}^{N \times N}$ be a matrix with $p \times p$ blocks, structured in the sense of \S \ref{sec:preliminaries:matrix-structure}, and $\epsilon > 0$ a target relative precision. We first outline a one-level matrix compression scheme:
 \begin{remunerate}
  \item
   For $i = 1, \dots, p$, use the ID to compress the row space of the $i$th off-diagonal block row to precision $\epsilon$. Let $L_{i}$ denote the corresponding row projection matrix.
  \item
   Similarly, for $j = 1, \dots, p$, use the ID to compress the column space of the $j$th off-diagonal block column to precision $\epsilon$. Let $R_{j}$ denote the corresponding column projection matrix.
  \item
   Approximate the off-diagonal blocks of $A$ by $A_{ij} \approx L_{i} S_{ij} R_{j}$ for $i \neq j$, where $S_{ij}$ is the submatrix of $A_{ij}$ defined by the row and column skeletons associated with $L_{i}$ and $R_{j}$, respectively.
 \end{remunerate}

\begin{figure}
  \begin{center}
   \includegraphics{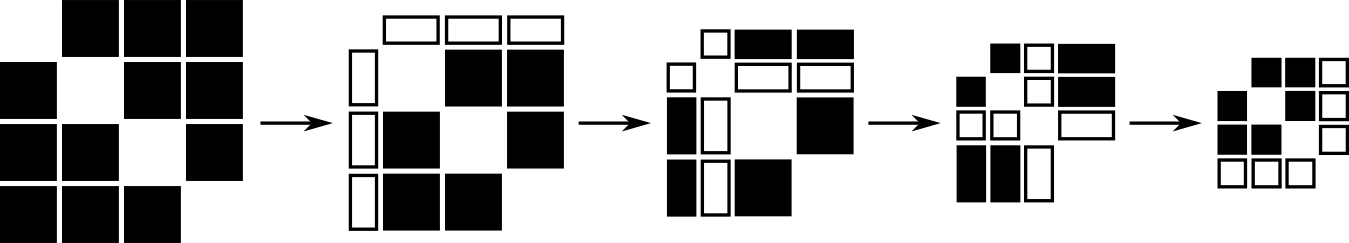}
  \end{center}
  \caption{One level of matrix compression, obtained by sequentially compressing the off-diagonal block rows and columns. At each step, the matrix blocks whose row or column spaces are being compressed are highlighted in white.}
  \label{fig:mat-comp}
 \end{figure}

 This yields precisely the matrix structure discussed in \S \ref{sec:preliminaries}, following (\ref{eq:block separable}). The one-level scheme is illustrated graphically in Fig.\ \ref{fig:mat-comp}.

\begin{figure}
  \begin{center}
   \includegraphics{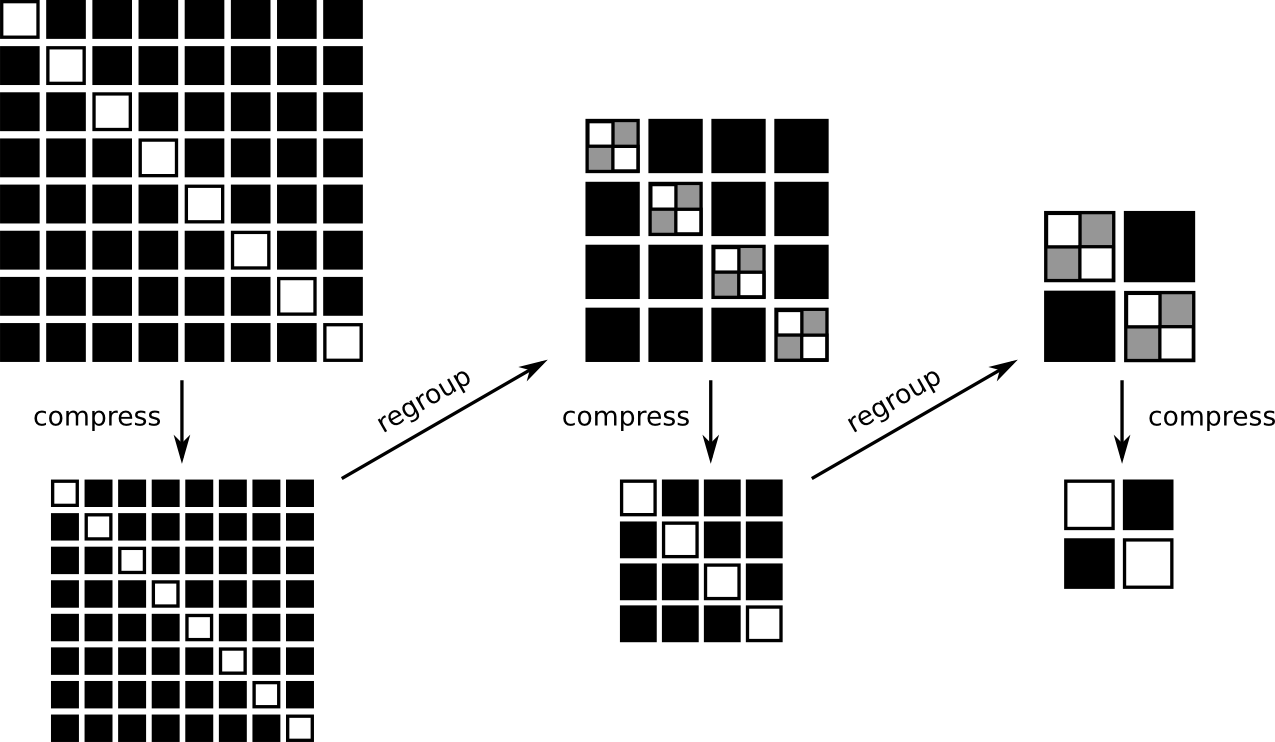}
  \end{center}
  \caption{Multilevel matrix compression, comprising alternating steps of compression and regrouping via ascension of the index tree. The diagonal blocks (white and gray) are not compressed, but are instead extracted at each level of the tree; they are shown here only to illustrate the regrouping process.}
  \label{fig:multi-comp}
 \end{figure}

 The multilevel algorithm is now just a simple extension based on the observation that by ascending one level in the index tree and regrouping blocks accordingly, we can compress the skeleton matrix $S$ in (\ref{eq:compressed-representation}) in exactly the same form, leading to a procedure that we naturally call {\em recursive skeletonization} (Fig.\ \ref{fig:multi-comp}).

 The full algorithm may be specified as follows:
 \begin{remunerate}
  \item
   Starting at the leaves of the tree, extract the diagonal blocks and perform one level of compression of the off-diagonal blocks.
  \item
   Move up one level in the tree and regroup the matrix blocks according to the tree structure. Terminate if the new level is the root; otherwise, extract the diagonal blocks, recompress the off-diagonal blocks, and repeat.
 \end{remunerate}

 The result is a telescoping representation
 \begin{equation}
  A \approx D^{\left( 1 \right)} + L^{\left( 1 \right)} \left[ D^{\left( 2 \right)} + L^{\left( 2 \right)} \left( \cdots D^{\left( \lambda \right)} + L^{\left( \lambda \right)} S R^{\left( \lambda \right)} \cdots \right) R^{\left( 2 \right)} \right] R^{\left( 1 \right)},
  \label{eq:multilevel-compression}
 \end{equation}
 where the superscript indexes the compression level $l = 1, \dots, \lambda$.

\begin{figure}
  \begin{center}
    \includegraphics{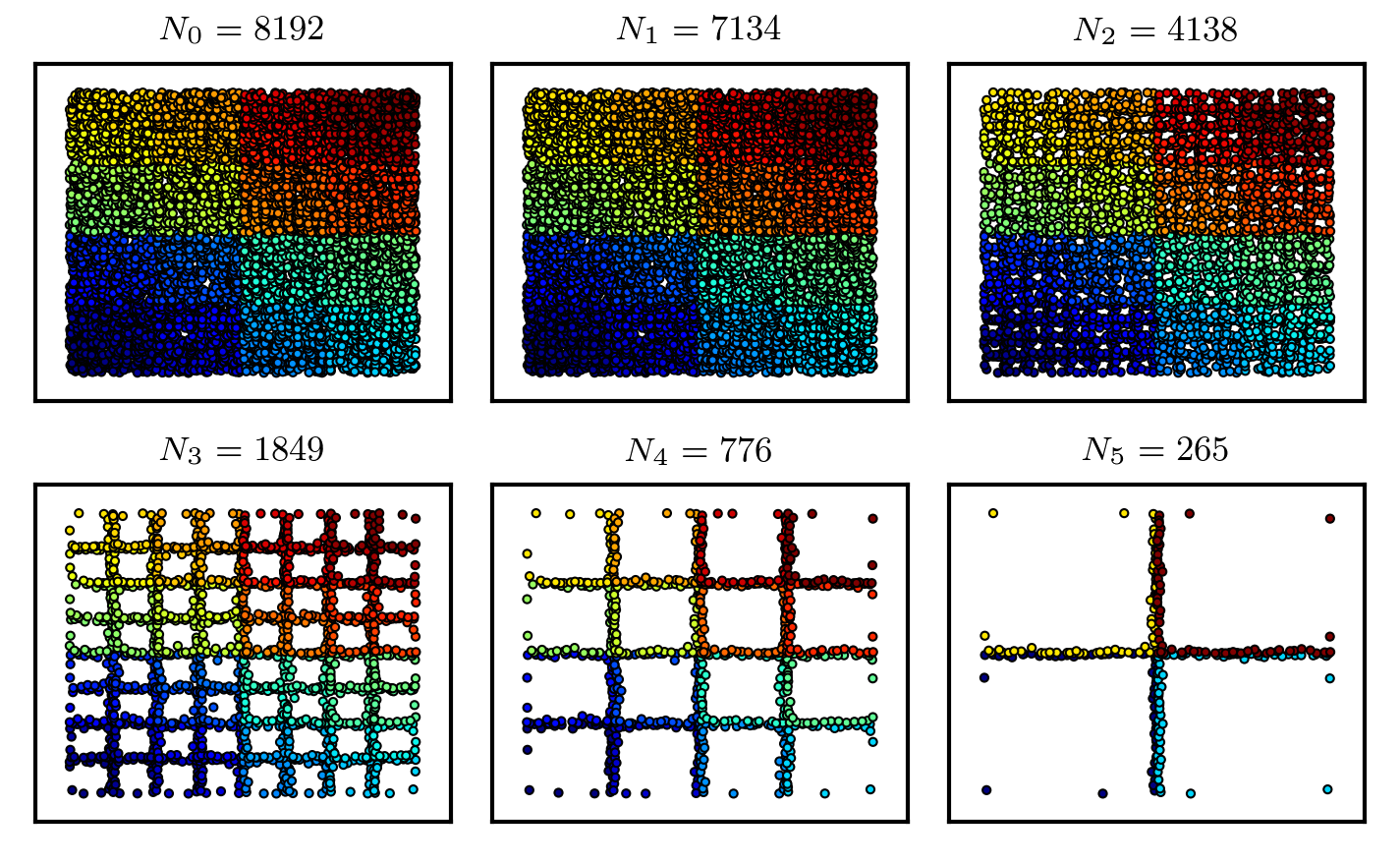}
  \end{center}
  \caption{Sparsification by recursive skeletonization. Logarithmic interactions between $N = 8192$ points in the unit square are compressed to relative precision $\epsilon = 10^{-3}$ using a five-level quadtree-based scheme. At each level, the surviving skeletons are shown, colored by block index, with the total number of skeletons remaining given by $N_{l}$ for compression level $l = 0, \dots, 5$, where $l = 0$ denotes the original uncompressed system.}
  \label{fig:sparsify}
 \end{figure}

 {\it Example 3.} As a demonstration of the multilevel compression technique, consider the matrix defined by $N = 8192$ points uniformly distributed in the unit square, interacting via the 2D Laplace Green's function (\ref{eq:2d-laplace}) and sorted according to a quadtree ordering. The sequence of skeletons remaining after each level of compression to $\epsilon = 10^{-3}$ is shown in Fig.\ \ref{fig:sparsify}, from which we see that compression creates a {\em sparsification} of the sources which, in a geometric setting, leaves skeletons along the boundaries of each block.

 \subsection{Accelerated compression via proxy surfaces}
 The computational cost of the algorithm described in the previous section is dominated by the fact that each step is global: that is, compressing the row or column space for each block requires accessing all other blocks in that row or column. If no further knowledge of the matrix is available, this is indeed necessary. However, as noted by \cite{cheng:2005:siam-j-sci-comput,greengard:2009:acta-numer,martinsson:2005:j-comput-phys,martinsson:2007:j-comput-phys}, this global work can often be replaced by a local one, resulting in considerable savings.

\begin{figure}
  \begin{center}
   \includegraphics{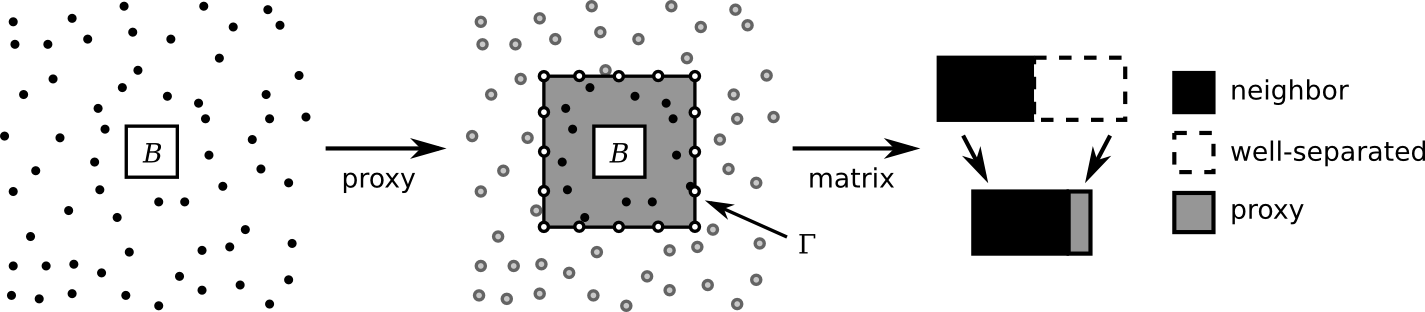}
  \end{center}
  \caption{Accelerated compression using proxy surfaces. The field within a region $B$ due to a distribution of exterior sources (left) can be decomposed into neighboring and well-separated contributions. By representing the latter via a proxy surface $\Gamma$ (center), the matrix dimension to compress against for the block corresponding to $B$ (right) can be reduced to the number of neighboring points plus a constant set of points on $\Gamma$, regardless of how many points lie beyond $\Gamma$.}
  \label{fig:proxy-comp}
 \end{figure}

 A sufficient condition for this acceleration is that the matrix correspond to evaluating a potential field for which some form of Green's identities hold. It is easiest to present the main ideas in the context of Laplace's equation. For this, consider Fig.\ \ref{fig:proxy-comp}, which depicts a set of sources in the plane. We assume that block index $i$ corresponds to the sources in the central square $B$. The $i$th off-diagonal block row then corresponds to the interactions of all points outside $B$ with all points inside $B$. We can separate this into contributions from the near neighbors of $B$, which are local, and the distant sources, which lie outside the near-neighbor domain, whose boundary is denoted by $\Gamma$. But any field induced by the distant sources induces a harmonic function inside $\Gamma$ and can therefore be replicated by a charge density on $\Gamma$ itself. Thus, rather than using the detailed structure of the distant points, the row (incoming) skeletons for $B$ can be extracted by considering just the combination of the near-neighbor sources and an artifical set of charges placed on $\Gamma$, which we refer to as a {\em proxy surface}. Likewise, the column (outgoing) skeletons for $B$ can be determined by considering only the near neighbors and the proxy surface. If the potential field is correct on the proxy surface, it will be correct at all more distant points (again via some variant of Green's theorem).

 The interaction rank between $\Gamma$ and $B$ is constant (depending on the desired precision) from standard multipole estimates \cite{greengard:1987:j-comput-phys,greengard:1997:acta-numer}. In summary, the number of points required to discretize $\Gamma$ is constant, and the dimension of the matrix to compress against for the block corresponding to $B$ is essentially just the number of points in the physically neighboring blocks.

 Similar arguments hold for other kernels of potential theory including the heat, Helmholtz, Yukawa, Stokes, and elasticity kernels, though care must be taken for oscillatory problems which could require a combination of single and double layer potentials to avoid spurious resonances in the representation for the exterior.

 \subsection{Compressed matrix-vector multiplication}
 The compressed representation (\ref{eq:multilevel-compression}) admits an obvious fast algorithm for computing the matrix-vector product $\mathbf{y} = A \mathbf{x}$.
As shown in \cite{gillman:2012:front-math-china},
one simply applies the matrices in (\ref{eq:multilevel-compression}) from right to left.
Like the FMM, this procedure can be thought of as occurring in two passes:
 \begin{remunerate}
  \item
   An upward pass, corresponding to the sequential application of the column projection matrices $R^{(l)}$, which hierarchically compress the input data $\mathbf{x}$ to the column (outgoing) skeleton subspace.
  \item
   A downward pass, corresponding to the sequential application of the row projection matrices $L^{(l)}$, which hierarchically project onto the row (incoming) skeleton subspace and, from there, back onto the output elements $\mathbf{y}$.
 \end{remunerate}

 \subsection{Compressed matrix inversion}
 \label{sec:inversion}
 The representation (\ref{eq:multilevel-compression}) also permits a fast algorithm for the direct inversion of nonsingular matrices. The one-level scheme was discussed in \S \ref{sec:preliminaries}. In the multilevel scheme, the system $S \mathbf{z} = \mathbf{y}$ in (\ref{eq:sparse-embedding}) is itself expanded in the same form, leading to the sparse embedding
 \begin{equation}
  \left[
  \begin{array}{ccccccc}
   D^{\left( 1 \right)} & L^{\left( 1 \right)}\\
   R^{\left( 1 \right)} & & -I\\
   & -I & D^{\left( 2 \right)} & L^{\left( 2 \right)}\\
   & & R^{\left( 2 \right)} & \ddots & \ddots\\
   & & & \ddots & D^{\left( \lambda \right)} & L^{\left( \lambda \right)}\\
   & & & & R^{\left( \lambda \right)} & & -I\\
   & & & & & -I & S
  \end{array} \right] \left[
  \begin{array}{c}
   \mathbf{x}\\
   \mathbf{y}^{\left( 1 \right)}\\
   \mathbf{z}^{\left( 1 \right)}\\
   \vdots\\
   \vdots\\
   \mathbf{y}^{\left( \lambda \right)}\\
   \mathbf{z}^{\left( \lambda \right)}
  \end{array} \right] = \left[
  \begin{array}{c}
   \mathbf{b}\\
   \mathbf{0}\\
   \mathbf{0}\\
   \vdots\\
   \vdots\\
   \mathbf{0}\\
   \mathbf{0}
  \end{array} \right].
  \label{eq:multilevel-embedding}
 \end{equation}

 To understand the consequences of this sparse representation, it is instructive to consider the special case in which the row and column skeleton dimensions are identical for each block, say $k$, so that the total row and column skeleton dimensions are $K \equiv K_{\row} = K_{\col} = pk$. Then, studying (\ref{eq:sparse-embedding}) first and assuming that $D$ is invertible, block elimination of $\mathbf{x}$ and $\mathbf{y}$ yields
 $$
  \left( \Lambda + S \right) \mathbf{z} = \Lambda R D^{-1} \mathbf{b},
 $$
 where $\Lambda = (R D^{-1} L)^{-1} \in \mathbb{C}^{K \times K}$ is block diagonal. Back substitution then yields
 $$
  \mathbf{x} = \left[ D^{-1} - D^{-1} L \Lambda R D^{-1} + D^{-1} L \Lambda \left( \Lambda + S \right)^{-1} \Lambda R D^{-1} \right] \mathbf{b}.
 $$
 In other words, the matrix inverse is
 \begin{equation}
  A^{-1} \approx \mathcal{D} + \mathcal{L} \mathcal{S}^{-1} \mathcal{R},
  \label{eq:compressed-inverse}
 \end{equation}
 where
 $$
  \mathcal{D} = D^{-1} - D^{-1} L \Lambda R D^{-1} \in \mathbb{C}^{N \times N}
 $$
 and
 $$
  \mathcal{L} = D^{-1} L \Lambda \in \mathbb{C}^{N \times K}, \qquad \mathcal{R} = \Lambda R D^{-1} \in \mathbb{C}^{K \times N}
 $$
 are all block diagonal, and
 $$
  \mathcal{S} = \Lambda + S \in \mathbb{C}^{K \times K}
 $$
 is dense. Note that $\mathcal{S}$ is equal to the skeleton matrix $S$ with its diagonal blocks filled in. Thus, (\ref{eq:compressed-inverse}) is a compressed representation of $A^{-1}$ with minimal fill-in over the original compressed representation (\ref{eq:compressed-representation}) of $A$. In the multilevel setting, one carries out the above factorization recursively, since $\mathcal{S}$ can now be inverted in the same manner:
 \begin{equation}
  A^{-1} \approx \mathcal{D}^{\left( 1 \right)} + \mathcal{L}^{\left( 1 \right)} \left[ \mathcal{D}^{\left( 2 \right)} + \mathcal{L}^{\left( 2 \right)} \left( \cdots \mathcal{D}^{\left( \lambda \right)} + \mathcal{L}^{\left( \lambda \right)} \mathcal{S}^{-1} \mathcal{R}^{\left( \lambda \right)} \cdots \right) \mathcal{R}^{\left( 2 \right)} \right] \mathcal{R}^{\left( 1 \right)}.
  \label{eq:multilevel-inverse}
 \end{equation}
 This point of view is elaborated in \cite{gillman:2012:front-math-china}.

 In the general case, the preceding procedure may fail if $D$ happens to be singular and (more generally) may be numerically unstable if care is not taken to stabilize the block elimination scheme using some sort of pivoting. Thus, rather than using the ``hand-rolled'' Gaussian elimination scheme of \cite{gillman:2011:univ-colorado-boulder,gillman:2012:front-math-china,martinsson:2005:j-comput-phys}
to compute the telescoping inverse (\ref{eq:multilevel-inverse}), we rely instead on the existence of high-quality sparse direct solver software. More precisely, we simply supply UMFPACK with the sparse representation (\ref{eq:multilevel-embedding}) and let it compute the corresponding factorization. Numerical results show that the performance is similar to that expected from (\ref{eq:multilevel-inverse}).

 \section{Complexity analysis}
 \label{sec:complexity-analysis}
 For the sake of completeness, we briefly analyze the complexity of the algorithm presented in \S \ref{sec:algorithm} for a typical example: discretization of the integral operator (\ref{eq:integral-operator}), where the integral kernel has smoothness properties similar to that of the Green's function for the Laplace equation. We follow the analysis of \cite{gillman:2011:univ-colorado-boulder,gillman:2012:front-math-china,martinsson:2005:j-comput-phys} and estimate costs for the ``hand-rolled'' Gaussian elimination scheme. We ignore quadrature issues and assume that we are given a matrix $A$ acting on $N$ points distributed randomly in a $d$-dimensional domain, sorted by an orthtree that uniformly subdivides until all block sizes are $\bigO (1)$. (In 1D, an orthtree is a binary tree; in 2D, it is a quadtree; and in 3D, it is an octree.)

 For each compression level $l = 1, \dots, \lambda$, with $l = 1$ being the finest, let $p_{l}$ be the number of matrix blocks, and $n_{l}$ and $k_{l}$ the uncompressed and compressed block dimensions, respectively, assumed equal for all blocks and identical across rows and columns, for simplicity. We first make the following observations:
 \begin{romannum}
  \item
   The total matrix dimension is $p_{1} n_{1} = N$, where $n_{1} = \bigO (1)$, so $p_{1} \sim N$.
  \item
   Each subdivision increases the number of blocks by a factor of roughly $2^{d}$, so $p_{l} \sim p_{l - 1} / 2^{d} \sim p_{1} / 2^{d(l - 1)}$. In particular, $p_{\lambda} = \bigO (1)$, so $\lambda \sim \log_{2^{d}} N = (1/d) \log N$.
  \item
   The total number of points at level $l > 1$ is equal to the total number of skeletons at level $l - 1$, i.e., $p_{l} n_{l} = p_{l - 1} k_{l - 1}$, so $n_{l} \sim 2^{d} k_{l - 1}$.
 \end{romannum}

\begin{figure}
  \begin{center}
   \includegraphics{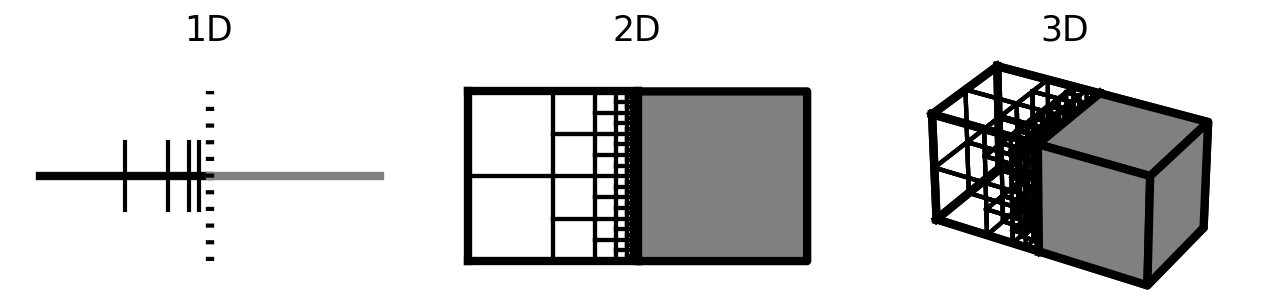}
  \end{center}
  \caption{The interaction rank between two adjacent blocks can be calculated by recursively subdividing the source block (white) into well-separated subblocks with respect to the target (gray), each of which have constant rank.}
  \label{fig:recur-subdiv}
 \end{figure}

 Furthermore, we note that $k_{l}$ is on the order of the interaction rank between two adjacent blocks at level $l$, which can be analyzed by recursive subdivision of the source block to expose well-separated structures with respect to the target (Fig.\ \ref{fig:recur-subdiv}).
 Assuming only that the interaction between a source subregion separated from a target by a distance of at least its own size is of constant rank (to a fixed precision $\epsilon$), we have
 $$
  k_{l} \sim \sum_{l = 1}^{\log_{2^{d}} n_{l}} 2^{\left( d - 1 \right) l} \sim \left\{
  \begin{array}{ll}
   \log n_{l} & \mbox{if $d = 1$},\\
   n_{l}^{1 - 1/d} & \mbox{if $d > 1$},
  \end{array} \right.
 $$
 where, clearly, $n_{l} \sim (p_{1} / p_{l}) n_{1} \sim 2^{d(l - 1)} n_{1}$, so
 $$
  k_{l} \sim \left\{
  \begin{array}{ll}
   \left( l - 1 \right) \log 2 + \log n_{1} & \mbox{if $d = 1$},\\
   2^{\left( d - 1 \right) \left( l - 1 \right)} n_{1}^{1 - 1/d} & \mbox{if $d > 1$}.
  \end{array} \right.
 $$

 \subsection{Matrix compression}
 From \S \ref{sec:preliminaries:interpolative-decomposition}, the cost of computing a rank-$k$ ID of an $m \times n$ matrix is $\bigO (mn \log k + k^{2} n)$. We will only consider the case of proxy compression, for which $m = \bigO (n_{l})$ for a block at level $l$, so the total cost is
 \begin{equation}
  T_{\cm} \sim \sum_{l = 1}^{\lambda} p_{l} \left( n_{l}^{2} \log k_{l} + k_{l}^{2} n_{l} \right) \sim \left\{
  \begin{array}{ll}
   N & \mbox{if $d = 1$},\\
   N^{3 \left( 1 - 1/d \right)} & \mbox{if $d > 1$}.
  \end{array} \right.
  \label{eq:complexity-cm}
 \end{equation}

 \subsection{Matrix-vector multiplication}
 \label{sec:complexity-analysis:mv}
 The cost of applying $D^{(l)}$ is $\bigO (p_{l} n_{l}^{2})$, while that of applying $L^{(l)}$ or $R^{(l)}$ is $\bigO (p_{l} k_{l} n_{l})$. Combined with the $\bigO ((p_{\lambda} k_{\lambda})^{2})$ cost of applying $S$, the total cost is
 \begin{equation}
  T_{\mv} \sim \sum_{l = 1}^{\lambda} p_{l} n_{l} \left( k_{l} + n_{l} \right) + \left( p_{\lambda} k_{\lambda} \right)^{2} \sim \left\{
  \begin{array}{ll}
   N & \mbox{if $d = 1$},\\
   N \log N & \mbox{if $d = 2$},\\
   N^{2 \left( 1 - 1/d \right)} & \mbox{if $d > 2$}.
  \end{array} \right.
  \label{eq:complexity-mv}
 \end{equation}

 \subsection{Matrix factorization and inverse application}
 \label{sec:complexity-analysis:sv}
 We turn now to the analysis of the cost of factorization using (\ref{eq:multilevel-inverse}). At each level $l$, the cost of constructing $D^{-1}$ and $\Lambda$ is $\bigO (p_{l} n_{l}^{3})$, after which forming $\mathcal{D}^{(l)}$, $\mathcal{L}^{(l)}$, and $\mathcal{R}^{(l)}$ all require $\bigO (p_{l} n_{l}^{2})$ operations; at the final level, the cost of constructing and inverting $\mathcal{S}$ is $\bigO ((p_{\lambda} k_{\lambda})^{3})$. Thus, the total cost is
 $$
  T_{\lu} \sim \sum_{l = 1}^{\lambda} p_{l} n_{l}^{3} + \left( p_{\lambda} k_{\lambda} \right)^{3},
 $$
 which has complexity (\ref{eq:complexity-cm}).

 Finally, we note that the dimensions of $\mathcal{D}^{(l)}$, $\mathcal{L}^{(l)}$, $\mathcal{R}^{(l)}$, and $\mathcal{S}^{-1}$ are the same as those of $D^{(l)}$, $L^{(l)}$, $R^{(l)}$, and $S$, respectively. Thus, the total cost of applying the inverse, denoted by $T_{\sv}$, has the same complexity as $T_{\mv}$, namely (\ref{eq:complexity-mv}).

 In our UMFPACK-based approach, the estimation of cost is a rather complicated task, and we do not attempt to carry out a detailed analysis of its performance. Suffice it to say, there is a one-to-one correspondence between the ``hand-rolled'' Gaussian elimination approach and one possible elimination scheme in UMFPACK. Since that solver is highly optimized, the asymptotic cost should be the same (or better). For some matrices, it is possible that straight Gaussian elimination may be unstable without pivoting, while UMFPACK will carry out a backward-stable scheme. This is a distinct advantage of the sparse matrix approach although the complexity and fill-in analysis then becomes more involved.

 \subsection{Storage}
 An important issue in direct solvers, of course, is that of storage requirements. In the present setting the relevant matrices are the compressed sparse representation (\ref{eq:multilevel-embedding}) and the factorization computed within UMFPACK. This will be (\ref{eq:complexity-mv}) for the forward operator and, in the absence of pivoting, for the sparse factorization as well. If pivoting is required, the analysis is more complex as it involves some matrix fill-in and is postponed to future work.

 \section{Error analysis}
 \label{sec:error-analysis}
 We now state some simple error estimates for applying and inverting a compressed matrix. Let $A$ be the original matrix and $A_{\epsilon}$ its compressed representation, constructed using the algorithm of \S \ref{sec:algorithm} such that
 $$
  \frac{\left\| A - A_{\epsilon} \right\|}{\left\| A \right\|} \leq \epsilon
 $$
 for some $\epsilon > 0$. Note that this need not be the same as the specified local precision in the ID since errors may accumulate across levels. However, as in \cite{gillman:2012:front-math-china}, we have found that such error propagation is mild.

 Let $\mathbf{x}$ and $\mathbf{b}$ be vectors such that $A \mathbf{x} = \mathbf{b}$. Then it is straightforward to verify that for $\mathbf{b}_{\epsilon} = A_\epsilon \mathbf{x}$,
 $$
  \frac{\left\| \mathbf{b} - \mathbf{b}_{\epsilon} \right\|}{\left\| \mathbf{b} \right\|} \leq \epsilon \left\| A \right\| \left\| A^{-1} \right\| = \epsilon \kappa \left( A \right),
 $$
 where $\kappa (A)$ is the condition number of $A$. Furthermore, if $\mathbf{x}_{\epsilon} = A_{\epsilon}^{-1} \mathbf{b}$, then
 $$
  \frac{\left\| \mathbf{x} - \mathbf{x}_{\epsilon} \right\|}{\left\| \mathbf{x} \right\|} \leq \frac{2 \epsilon \kappa \left( A \right)}{1 - \epsilon \kappa \left( A \right)}.
 $$
 In particular, if $A$ is well-conditioned, e.g., $A$ is the discretization of a second-kind integral equation, then $\kappa (A) = \bigO (1)$, so
 $$
  \frac{\left\| \mathbf{x} - \mathbf{x}_{\epsilon} \right\|}{\left\| \mathbf{x} \right\|} = \bigO \left( \epsilon \right).
 $$

 \section{Numerical examples}
 \label{sec:numerical-examples}
 In this section, we investigate the efficiency and flexibility of our algorithm by considering some representative examples. We begin with timing benchmarks for the Laplace and Helmholtz kernels in 2D and 3D, using the algorithm both as an FMM and as a direct solver, followed by applications in molecular electrostatics and multiple scattering.

 All matrices were blocked using quadtrees in 2D and octrees in 3D, uniformly subdivided until all block sizes were $\bigO (1)$, but adaptively truncating empty boxes during the refinement process. Only proxy compression was considered, with proxy surfaces constructed on the boundary of the supercell enclosing the neighbors of each block. We discretized all proxy surfaces using a constant number of points, independent of the matrix size $N$: for the Laplace equation, this constant depended only on the compression precision $\epsilon$, while for the Helmholtz equation, it depended also on the wave frequency, chosen to be consistent with the Nyquist-Shannon sampling theorem. Computations were performed over $\mathbb{R}$ instead of $\mathbb{C}$, where possible. The algorithm was implemented in Fortran, and all experiments were performed on a 2.66 GHz processor in double precision.

 In many instances, we compare our results against those obtained using LAPACK/ATLAS \cite{anderson:1999:soc-indust-appl-math,whaley:2001:parallel-comput} and the FMM \cite{cheng:2006:j-comput-phys,greengard:1987:j-comput-phys,greengard:1997:acta-numer,rokhlin:1990:j-comput-phys}. All FMM timings were computed using the open-source FMMLIB package \cite{gimbutas:in-prep}, which is a fairly efficient implementation but does not include the plane-wave optimizations of \cite{cheng:2006:j-comput-phys,greengard:1997:acta-numer} or the diagonal translation operators of \cite{rokhlin:1990:j-comput-phys}.

 \subsection{Generalized fast multipole method}
 \label{sec:numerical-examples:generalized-fmm}
 We first consider the use of recursive skeletonization as a generalized FMM for the rapid computation of matrix-vector products.

 \subsubsection{The Laplace equation}
 We considered two point distributions in the plane: points on the unit circle and in the unit square, hereafter referred to as the 2D surface and volume cases, respectively. We assumed that the governing matrix corresponds to the interaction of charges via the Green's function (\ref{eq:2d-laplace}). The surface case is typical of layer-potential evaluation when using boundary integral equations. Since a domain boundary in 2D can be described by a single parameter (such as arclength), it is a 1D domain, so the expected complexities from \S \ref{sec:complexity-analysis} correspond to $d = 1$: $\bigO (N)$ work for both matrix compression and application. (See \cite{gillman:2012:front-math-china} for a detailed discussion of the $d = 1$ case.) In the volume case, the dimension is $d = 2$, so the expected complexities are $\bigO (N^{3/2})$ and $\bigO (N \log N)$ for compression and application, respectively.

 For the 3D Laplace kernel (\ref{eq:3d-laplace}), we considered surface and volume point geometries on the unit sphere and within the unit cube, respectively. The corresponding dimensions are $d = 2$ and $d = 3$. Thus, the expected complexities for the 3D surface case are $\bigO (N^{3/2})$ and $\bigO (N \log N)$ for compression and application, respectively, while those for the 3D volume case are $\bigO (N^{2})$ and $\bigO (N^{4/3})$, respectively.

\begin{figure}
  \begin{center}
   \includegraphics{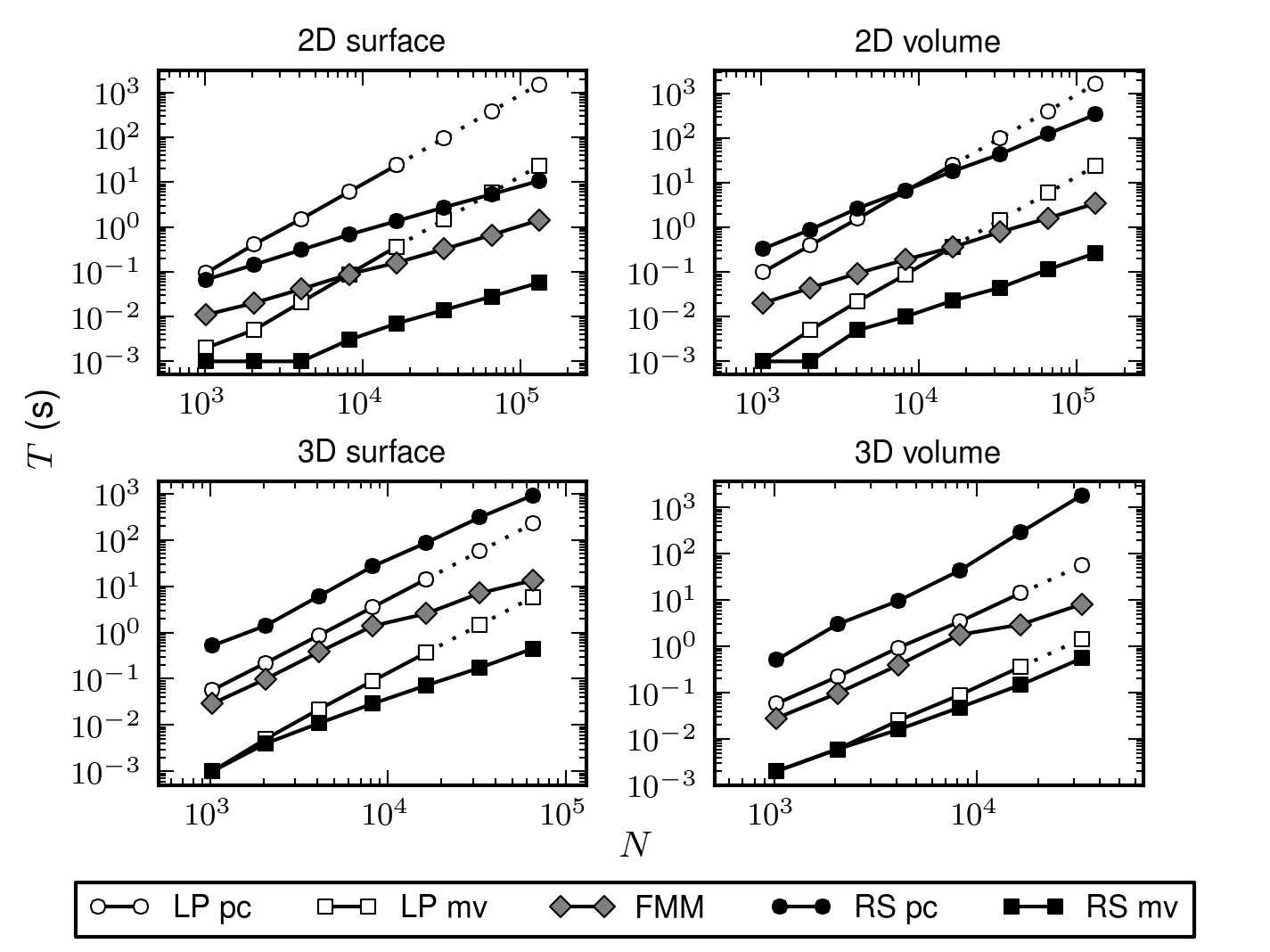}
  \end{center}
  \caption{CPU times for applying the Laplace kernel in various cases using LAPACK/ATLAS (LP), the FMM, and recursive skeletonization (RS) as a function of the matrix size $N$. For LP and RS, the computation is split into two parts: precomputation (pc), for LP consisting of matrix formation and for RS of matrix compression, and matrix-vector multiplication (mv). The precision of the FMM and RS was set at $\epsilon = 10^{-9}$. Dotted lines indicate extrapolated data.}
  \label{fig:apply-lap}
 \end{figure}

 We present timing results for each case and compare with LAPACK/ATLAS and the FMM for a range of $N$ at $\epsilon = 10^{-9}$. Detailed data are provided in Tables \ref{tab:apply-lap-2ds}--\ref{tab:apply-lap-3dv} and plotted in Fig.\ \ref{fig:apply-lap}.
 It is evident that our algorithm scales as predicted.
 \begin{table}
  \caption{Numerical results for applying the Laplace kernel in the 2D surface case at precision $\epsilon = 10^{-9}$: $N$, uncompressed matrix dimension; $K_{\row}$, row skeleton dimension; $K_{\col}$, column skeleton dimension; $T_{\cm}$, matrix compression time (s); $T_{\mv}$, matrix-vector multiplication time (s); $E$, relative error; $M$, required storage for compressed matrix (MB).}
  \label{tab:apply-lap-2ds}
  \begin{center}
   {\footnotesize
    \begin{tabular}{r|rrcccc}
     \hline
     \multicolumn{1}{c|}{$N$} & \multicolumn{1}{c}{$K_{\row}$} & \multicolumn{1}{c}{$K_{\col}$} & $T_{\cm}$ & $T_{\mv}$ & $E$ & $M$\\
     \hline
       1024 &  94 &  94 & 6.7\E$-$2 & 1.0\E$-$3 & 3.1\E$-$8 & 8.5\E$-$1\\
       2048 & 105 & 104 & 1.4\E$-$1 & 1.0\E$-$3 & 4.5\E$-$8 & 1.7\E$+$0\\
       4096 & 113 & 114 & 3.1\E$-$1 & 1.0\E$-$3 & 1.1\E$-$7 & 3.4\E$+$0\\
       8192 & 123 & 123 & 6.7\E$-$1 & 3.0\E$-$3 & 4.4\E$-$7 & 6.4\E$+$0\\
      16384 & 133 & 134 & 1.4\E$+$0 & 7.0\E$-$3 & 4.0\E$-$7 & 1.3\E$+$1\\
      32768 & 142 & 142 & 2.7\E$+$0 & 1.4\E$-$2 & 4.7\E$-$7 & 2.5\E$+$1\\
      65536 & 150 & 149 & 5.4\E$+$0 & 2.8\E$-$2 & 9.4\E$-$7 & 5.0\E$+$1\\
     131072 & 159 & 158 & 1.1\E$+$1 & 5.7\E$-$2 & 9.8\E$-$7 & 1.0\E$+$2\\
     \hline
    \end{tabular}
   }
  \end{center}
 \end{table}
 Its performance in 2D is particularly strong.
 \begin{table}
  \caption{Numerical results for applying the Laplace kernel in the 2D volume case at precision $\epsilon = 10^{-9}$; notation as in Table \ref{tab:apply-lap-2ds}.}
  \label{tab:apply-lap-2dv}
  \begin{center}
   {\footnotesize
    \begin{tabular}{r|rrcccc}
     \hline
     \multicolumn{1}{c|}{$N$} & \multicolumn{1}{c}{$K_{\row}$} & \multicolumn{1}{c}{$K_{\col}$} & $T_{\cm}$ & $T_{\mv}$ & $E$ & $M$\\
     \hline
       1024 &  299 &  298 & 3.3\E$-$1 & 1.0\E$-$3 & 3.6\E$-$10 & 2.9\E$+$0\\
       2048 &  403 &  405 & 8.9\E$-$1 & 1.0\E$-$3 & 3.7\E$-$10 & 7.1\E$+$0\\
       4096 &  570 &  570 & 2.7\E$+$0 & 5.0\E$-$3 & 1.0\E$-$09 & 1.8\E$+$1\\
       8192 &  795 &  793 & 6.8\E$+$0 & 1.0\E$-$2 & 8.8\E$-$10 & 4.3\E$+$1\\
      16384 & 1092 & 1091 & 1.8\E$+$1 & 2.3\E$-$2 & 7.7\E$-$10 & 1.0\E$+$2\\
      32768 & 1506 & 1505 & 4.4\E$+$1 & 4.5\E$-$2 & 1.0\E$-$09 & 2.3\E$+$2\\
      65536 & 2099 & 2101 & 1.3\E$+$2 & 1.1\E$-$1 & 1.1\E$-$09 & 5.3\E$+$2\\
     131072 & 2904 & 2903 & 3.4\E$+$2 & 2.7\E$-$1 & 1.1\E$-$09 & 1.2\E$+$3\\
     \hline
    \end{tabular}
   }
  \end{center}
 \end{table}
 Not only does our algorithm beat the $\bigO (N^{2})$ uncompressed matrix-vector product for modest $N$, it is faster even than the $\bigO (N)$ FMM (at least after compression).
 \begin{table}
  \caption{Numerical results for applying the Laplace kernel in the 3D surface case at precision $\epsilon = 10^{-9}$; notation as in Table \ref{tab:apply-lap-2ds}.}
  \label{tab:apply-lap-3ds}
  \begin{center}
   {\footnotesize
    \begin{tabular}{r|rrcccc}
     \hline
     \multicolumn{1}{c|}{$N$} & \multicolumn{1}{c}{$K_{\row}$} & \multicolumn{1}{c}{$K_{\col}$} & $T_{\cm}$ & $T_{\mv}$ & $E$ & $M$\\
     \hline
      1024 &  967 &  967 & 5.2\E$-$1 & 1.0\E$-$3 & 1.0\E$-$11 & 7.7\E$+$0\\
      2048 & 1531 & 1532 & 1.4\E$+$0 & 4.0\E$-$3 & 1.8\E$-$10 & 2.2\E$+$1\\
      4096 & 2298 & 2295 & 6.1\E$+$0 & 1.1\E$-$2 & 1.4\E$-$10 & 6.2\E$+$1\\
      8192 & 3438 & 3426 & 2.7\E$+$1 & 2.9\E$-$2 & 1.2\E$-$10 & 1.7\E$+$2\\
     16384 & 4962 & 4950 & 8.7\E$+$1 & 7.2\E$-$2 & 3.0\E$-$10 & 4.2\E$+$2\\
     32768 & 6974 & 6987 & 3.1\E$+$2 & 1.7\E$-$1 & 4.3\E$-$10 & 9.9\E$+$2\\
     65536 & 9899 & 9925 & 9.2\E$+$2 & 4.5\E$-$1 & 7.7\E$-$10 & 2.3\E$+$3\\
     \hline
    \end{tabular}
   }
  \end{center}
 \end{table}
 In 3D, the same is true over the range of $N$ tested, although the increase in asymptotic complexity would eventually make the scheme less competitive.
 \begin{table}
  \caption{Numerical results for applying the Laplace kernel in the 3D volume case at precision $\epsilon = 10^{-9}$; notation as in Table \ref{tab:apply-lap-2ds}.}
  \label{tab:apply-lap-3dv}
  \begin{center}
   {\footnotesize
    \begin{tabular}{r|rrcccc}
     \hline
     \multicolumn{1}{c|}{$N$} & \multicolumn{1}{c}{$K_{\row}$} & \multicolumn{1}{c}{$K_{\col}$} & $T_{\cm}$ & $T_{\mv}$ & $E$ & $M$\\
     \hline
      1024 &  1024 &  1024 & 5.1\E$-$1 & 2.0\E$-$3 & 9.3\E$-$16 & 8.4\E$+$0\\
      2048 &  1969 &  1969 & 3.0\E$+$0 & 6.0\E$-$3 & 5.6\E$-$12 & 3.2\E$+$1\\
      4096 &  3285 &  3287 & 9.7\E$+$0 & 1.6\E$-$2 & 6.8\E$-$11 & 9.8\E$+$1\\
      8192 &  5360 &  5362 & 4.4\E$+$1 & 4.8\E$-$2 & 6.3\E$-$11 & 3.0\E$+$2\\
     16384 &  8703 &  8707 & 2.9\E$+$2 & 1.5\E$-$1 & 5.7\E$-$11 & 9.3\E$+$2\\
     32768 & 14015 & 14013 & 1.9\E$+$3 & 5.5\E$-$1 & 7.5\E$-$11 & 2.9\E$+$3\\
     \hline
    \end{tabular}
   }
  \end{center}
 \end{table}
 In all cases studied, the compression time $T_{\cm}$ was larger than the time to apply the FMM by one (2D surface) to two (all other cases) orders of magnitude, while the compressed matrix-vector product time $T_{\mv}$ was consistently smaller by the same amount. Thus, our algorithm also shows promise as a fast iterative solver for problems requiring more than $\sim 10$--$100$ iterations. Furthermore, we note the effectiveness of compression: for $N = 131072$, the storage requirement for the uncompressed matrix is $137$ GB, whereas that for the compressed representations are only $100$ MB and $1.2$ GB in the 2D surface and volume cases, respectively; at a lower precision of $\epsilon = 10^{-3}$, these become just $40$ and $180$ MB. Finally, to provide some intuition about the behavior of the algorithm as a function of precision, we report the following timings for the 2D volume case with $N = 131072$: for $\epsilon = 10^{-3}$, $T_{\cm} = 41$ s and $T_{\mv} = 0.09$ s; for $\epsilon = 10^{-6}$, $T_{\cm} = 161$ s and $T_{\mv} = 0.18$ s; and for $\epsilon = 10^{-9}$, $T_{\cm} = 339$ s and $T_{\mv} = 0.27$ s.

 \subsubsection{The Helmholtz equation}
 We next considered the 2D and 3D Helmholtz kernels
 \begin{equation}
  G \left( x, y \right) = \frac{\imath}{4} H_{0}^{\left( 1 \right)} \left( k \left| x - y \right| \right)
  \label{eq:2d-helmholtz}
 \end{equation}
 and
 \begin{equation}
  G \left( x, y \right) = \frac{e^{\imath k \left| x - y \right|}}{4 \pi \left| x - y \right|},
  \label{eq:3d-helmholtz}
 \end{equation}
 respectively, where $H_{0}^{(1)}$ is the zeroth order Hankel function of the first kind and $k$ is the wavenumber. We used the same representative geometries as for the Laplace equation. The size of each domain $\Omega$ in wavelengths was given by
 $$
  \omega = \frac{k}{2 \pi} \diam(\Omega).
 $$

\begin{figure}
  \begin{center}
   \includegraphics{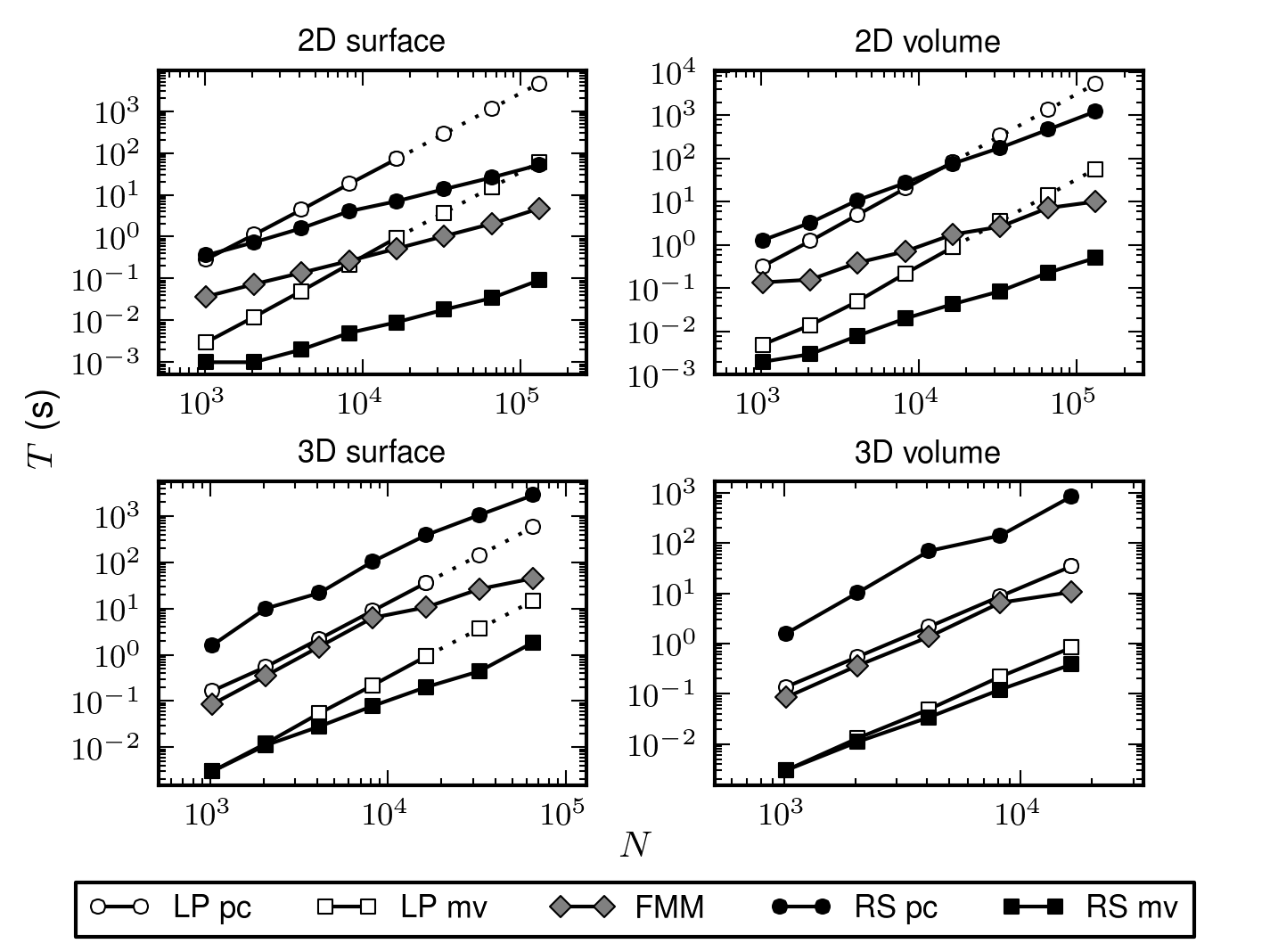}
  \end{center}
  \caption{CPU times for applying the Helmholtz kernel in various cases at low frequency ($\omega = 10$ in 2D and $\omega = 5$ in 3D) using LAPACK/ATLAS, the FMM, and recursive skeletonization at precision $\epsilon = 10^{-9}$; notation as in Fig.\ \ref{fig:apply-lap}.}
  \label{fig:apply-helm}
 \end{figure}

 Timing results against LAPACK/ATLAS and the FMM at low frequency ($\omega = 10$ in 2D and $\omega = 5$ in 3D) with $\epsilon = 10^{-9}$ are shown in Fig.\ \ref{fig:apply-helm}.
 In this regime, the performance is very similar to that for the Laplace equation, as both kernels are essentially non-oscillatory; detailed data are therefore omitted. However, as discussed in \cite{martinsson:2005:j-comput-phys}, the compression efficiency deteriorates as $\omega$ increases, due to the growing ranks of the matrix blocks. In the high-frequency regime, there is no asymptotic gain in efficiency. Still, numerical results suggest that the algorithm remains viable up to $\omega \sim 200$ in 2D and $\omega \sim 10$ in 3D. In all cases, the CPU times and storage requirements are larger than those for the Laplace equation by a factor of about two since all computations are performed over $\mathbb{C}$ instead of $\mathbb{R}$; in 2D, there is also the additional expense of computing $H_{0}^{(1)}$.

 \subsection{Recursive skeletonization as a direct solver}
 In this section, we study the behavior of our algorithm as a fast direct solver. More specifically, we considered the interior Dirichlet problem for the Laplace and Helmholtz equations in 2D and 3D, recast as a second-kind boundary integral equation using the double-layer representation (\ref{eq:double-layer}). Contour integrals in 2D were discretized using the trapezoidal rule, while surface integrals in 3D were discretized using Gaussian quadrature on flat triangles. In each case, we took as boundary data the field generated by an exterior point source; the error was assessed by comparing the field evaluated using the numerical solution via (\ref{eq:double-layer}) against the exact field due to that source at an interior location. As a benchmark, we also solved each system directly using LAPACK/ATLAS, as well as iteratively using GMRES with matrix-vector products accelerated by the FMM.

 \subsubsection{The Laplace equation}
\begin{figure}
  \begin{center}
   \includegraphics{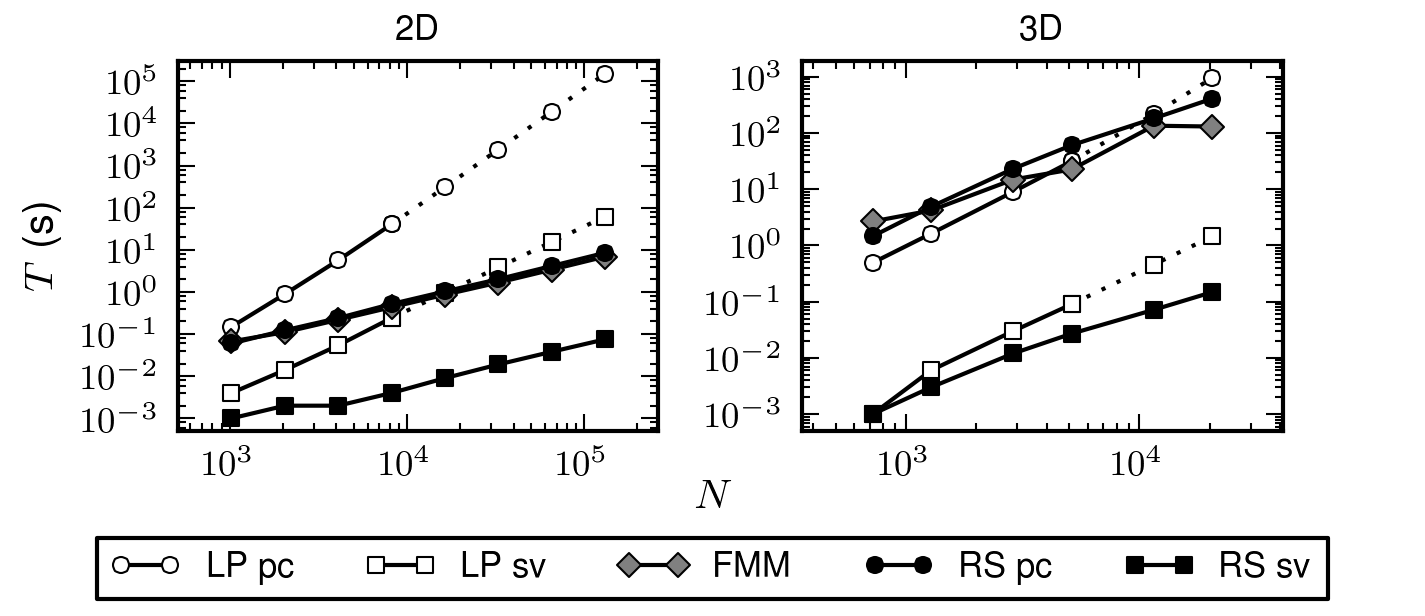}
  \end{center}
  \caption{CPU times for solving Laplace's equation in various cases using LAPACK/ATLAS (LP), FMM/GMRES (FMM), and recursive skeletonization (RS) as a function of the system size $N$. For LP and RS, the computation is split into two parts: precomputation (pc), for LP consisting of matrix formation and factorization, and for RS of matrix compression and factorization; and system solution (sv), consisting of matrix inverse application. The precision of the FMM and RS was set at $\epsilon = 10^{-9}$ in 2D and $\epsilon = 10^{-6}$ in 3D. Dotted lines indicate extrapolated data.}
  \label{fig:solve-lap}
 \end{figure}

 For the Laplace equation (\ref{eq:laplace-dirichlet}), the Green's function $G$ in (\ref{eq:2k-integral-eqn}) is given by (\ref{eq:2d-laplace}) in 2D and (\ref{eq:3d-laplace}) in 3D. As a model geometry, we considered an ellipse with aspect ratio $\alpha = 2$ (semi-major and -minor axes $a = 2$ and $b = 1$, respectively) in 2D and the unit sphere in 3D; these boundaries have dimensions $d = 1$ and $d = 2$, respectively. Timing results are shown in Fig.\ \ref{fig:solve-lap}, with detailed data given in Tables \ref{tab:solve-lap-2d} and \ref{tab:solve-lap-3d}; the precision was set to $\epsilon = 10^{-9}$ in 2D and $\epsilon = 10^{-6}$ in 3D.

 In 2D, the solver has linear complexity and is exceptionally fast, handily beating the $\bigO (N^{3})$ uncompressed direct solver, but also coming very close to the $\bigO (N)$ FMM/GMRES iterative solver. At $N = 131072$, for example, the total solution time for the recursive skeletonization algorithm was $T_{\RS} = 8.5$ s, while that for FMM/GMRES was $T_{\FMM} = 6.9$ s using $n_{\FMM} = 7$ iterations.
 \begin{table}
  \caption{Numerical results for solving Laplace's equation in 2D at precision $\epsilon = 10^{-9}$: $N$, uncompressed matrix dimension; $K_{\row}$, row skeleton dimension; $K_{\col}$, column skeleton dimension; $T_{\cm}$, matrix compression time (s); $T_{\lu}$, sparse matrix factorization time (s); $T_{\sv}$, inverse application time (s); $E$, relative error; $M$, required storage for compressed matrix inverse (MB).}
  \label{tab:solve-lap-2d}
  \begin{center}
   {\footnotesize
    \begin{tabular}{r|rrccccc}
     \hline
     \multicolumn{1}{c|}{$N$} & \multicolumn{1}{c}{$K_{\row}$} & \multicolumn{1}{c}{$K_{\col}$} & $T_{\cm}$ & $T_{\lu}$ & $T_{\sv}$ & $E$ & $M$\\
     \hline
       1024 & 30 & 30 & 3.4\E$-$2 & 2.5\E$-$2 & 1.0\E$-$3 & 9.0\E$-$11 & 1.6\E$+$0\\
       2048 & 29 & 30 & 7.0\E$-$2 & 5.1\E$-$2 & 2.0\E$-$3 & 9.0\E$-$12 & 3.3\E$+$0\\
       4096 & 30 & 30 & 1.4\E$-$1 & 9.8\E$-$2 & 2.0\E$-$3 & 8.3\E$-$11 & 6.8\E$+$0\\
       8192 & 30 & 31 & 3.0\E$-$1 & 2.1\E$-$1 & 4.0\E$-$3 & 1.6\E$-$10 & 1.4\E$+$1\\
      16384 & 31 & 31 & 5.5\E$-$1 & 4.5\E$-$1 & 9.0\E$-$3 & 5.5\E$-$10 & 2.8\E$+$1\\
      32768 & 30 & 30 & 1.1\E$+$0 & 8.5\E$-$1 & 1.9\E$-$2 & 4.9\E$-$12 & 5.6\E$+$1\\
      65536 & 30 & 30 & 2.3\E$+$0 & 1.8\E$+$0 & 3.8\E$-$2 & 1.1\E$-$11 & 1.1\E$+$2\\
     131072 & 29 & 29 & 4.6\E$+$0 & 3.7\E$+$0 & 7.5\E$-$2 & 8.5\E$-$11 & 2.2\E$+$2\\
     \hline
    \end{tabular}
   }
  \end{center}
 \end{table}
 It is worth emphasizing, however, that our solver is direct and possesses obvious advantages over FMM/GMRES, as described in \S \ref{sec:introduction}; in particular, the algorithm is relatively insensitive to geometric ill-conditioning.
 \begin{table}
  \caption{Numerical results for solving Laplace's equation in 3D at precision $\epsilon = 10^{-6}$; notation as in Table \ref{tab:solve-lap-2d}.}
  \label{tab:solve-lap-3d}
  \begin{center}
   {\footnotesize
    \begin{tabular}{r|rrccccc}
     \hline
     \multicolumn{1}{c|}{$N$} & \multicolumn{1}{c}{$K_{\row}$} & \multicolumn{1}{c}{$K_{\col}$} & $T_{\cm}$ & $T_{\lu}$ & $T_{\sv}$ & $E$ & $M$\\
     \hline
       720 &  628 &  669 & 1.3\E$+$0 & 1.1\E$-$1 & 1.0\E$-$3 & 9.8\E$-$5 & 4.6\E$+$0\\
      1280 &  890 &  913 & 4.5\E$+$0 & 4.0\E$-$1 & 3.0\E$-$3 & 5.5\E$-$5 & 1.1\E$+$1\\
      2880 & 1393 & 1400 & 2.1\E$+$1 & 2.0\E$+$0 & 1.2\E$-$2 & 2.4\E$-$5 & 5.5\E$+$1\\
      5120 & 1886 & 1850 & 5.5\E$+$1 & 5.4\E$+$0 & 2.7\E$-$2 & 1.3\E$-$5 & 1.3\E$+$2\\
     11520 & 2750 & 2754 & 1.6\E$+$2 & 1.7\E$+$1 & 7.2\E$-$2 & 6.2\E$-$6 & 3.5\E$+$2\\
     20480 & 3592 & 3551 & 3.7\E$+$2 & 4.1\E$+$1 & 1.5\E$-$1 & 3.3\E$-$6 & 6.9\E$+$2\\
     \hline
    \end{tabular}
   }
  \end{center}
 \end{table}
 Indeed, the direct solver edged out FMM/GMRES even at modest aspect ratios (for $N = 8192$ at $\epsilon = 10^{-12}$ with $\alpha = 8$: $T_{\RS} = 0.76$ s, $T_{\FMM} = 0.98$ s, $n_{\FMM} = 15$); for larger $\alpha$, the effect was even more pronounced ($\alpha = 512$: $T_{\RS} = 2.5$ s, $T_{\FMM} = 3.9$ s, $n_{\FMM} = 44$). Furthermore, the compressed inverse representation allows subsequent solves to be performed extremely rapidly; for instance, at $N = 131072$, the solve time was just $T_{\sv} = 0.07$ s, i.e., $T_{\FMM} / T_{\sv} \sim 100$. Thus, our algorithm is especially efficient in regimes where $T_{\sv}$ dominates (see, e.g., \cite{martinsson:2006:j-comput-phys}). Finally, we remark that although direct methods are traditionally very memory-intensive, our algorithm appears quite manageable in this regard: at $N = 131072$, the storage required for the compressed inverse was only $106$ MB for $\epsilon = 10^{-3}$, $172$ MB for $\epsilon = 10^{-6}$, and $222$ MB for $\epsilon = 10^{-9}$.

 In 3D, our solver has complexity $\bigO (N^{3/2})$. Hence, asymptotics dictate that it must eventually lose. However, our results demonstrate that even up to $N = 20480$, the solver remains surprisingly competitive. For example, at $N = 20480$, $T_{\RS} = 409$ s, while $T_{\FMM} = 131$ s with $n_{\FMM} = 3$; at $\epsilon = 10^{-9}$, the difference is almost negligible: $T_{\RS} = 850$ s, $T_{\FMM} = 839$ s, $n_{\FMM} = 5$. Thus, our algorithm remains a viable alternative for medium-scale problems. It is important to note that the {\em solve time} advantage is not lost even for large $N$, since the cost of each solve is only $\bigO (N \log N)$. In fact, the advantage is, remarkably, even more striking than in 2D: at $N = 20480$, $T_{\FMM} / T_{\sv} \sim 1000$; for $\epsilon = 10^{-9}$, $T_{\FMM} / T_{\sv} \sim 2500$.

 \subsubsection{The Helmholtz equation}
 We then considered the Helmholtz equation
 $$
  \left( \Delta + k^{2} \right) u = 0 \quad \mbox{in $\Omega$}, \qquad u = f \quad \mbox{on $\partial \Omega$},
 $$
 recast as a boundary integral equation (\ref{eq:2k-integral-eqn}), with Green's function (\ref{eq:2d-helmholtz}) in 2D and (\ref{eq:3d-helmholtz}) in 3D. This representation does not work for all frequencies, encountering spurious discrete resonances for $k$ beyond a critical value. We ignore that (well-understood) issue here and assume that the integral equation we obtain is invertible. The method itself does not falter in such cases, as discussed in \cite{martinsson:2005:j-comput-phys}.

 We used the same geometries and precisions as for the Laplace equation. In 2D, the double-layer kernel is weakly singular, so we modified the trapezoidal rule with tenth-order endpoint corrections \cite{kapur:1997:siam-j-numer-anal}. The frequency was set to $\omega = 10$ in 2D and $\omega = 3.18$ in 3D.

\begin{figure}
  \begin{center}
   \includegraphics{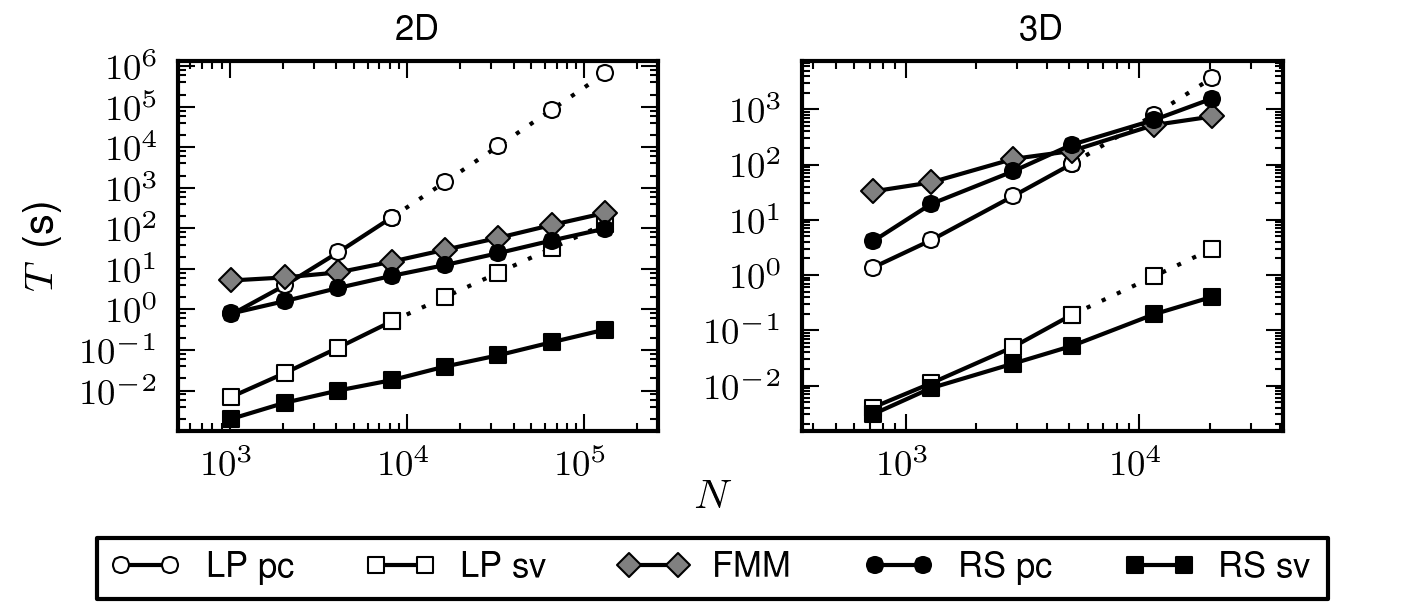}
  \end{center}
  \caption{CPU times for solving the Helmholtz equation in various cases at low frequency ($\omega = 10$ in 2D and $\omega = 3.18$ in 3D) using LAPACK/ATLAS, FMM/GMRES, and recursive skeletonization; notation as in Fig.\ \ref{fig:solve-lap}. The precision was set to $\epsilon = 10^{-9}$ in 2D and
$\epsilon = 10^{-6}$ in 3D.}
  \label{fig:solve-helm}
 \end{figure}

 Timing results are shown in Fig.\ \ref{fig:solve-helm}.
 The data are very similar to that for the Laplace equation, but with the direct solver actually beating FMM/GMRES in 2D. This is because the number of iterations required for FMM/GMRES scales as $n_{\FMM} = \bigO (\omega)$. Interestingly, even at moderately high frequencies, where we would expect the direct solver to break down as discussed in \S \ref{sec:numerical-examples:generalized-fmm}, the performance drop is more than compensated for by the increase in the number $n_{\FMM}$ of iterations. In short, we find that recursive skeletonization is faster than FMM/GMRES at low to moderate frequencies, provided that the memory requirement is not excessive.

 The story is much the same in 3D and the compressed solve time is again very fast: at $N = 20480$, $T_{\FMM} / T_{\sv} \sim 2000$.

 \subsection{Molecular electrostatics}
 An important application area for our solver is molecular electrostatics. A simplified model for this involves consideration of a molecular surface $\Sigma$, dividing $\mathbb{R}^{3}$ into $\Omega_{1}$ and $\Omega_{2}$, denoting the molecule and the solvent, respectively. We also suppose that the molecule has interior charges of strengths $q_{i}$ at locations $x_{i} \in \Omega_{1}$ for $i = 1, \dots, n$. The electrostatic potential $\varphi$ (ignoring salt effects in the solvent) then satisfies the Poisson equation
 $$
  -\nabla \cdot \left[ \varepsilon \left( x \right) \nabla \varphi \left( x \right) \right] = \sum_{i = 1}^{n} q_{i} \delta \left( x - x_{i} \right),
 $$
 where
 $$
  \varepsilon \left( x \right) = \left\{
  \begin{array}{ll}
   \varepsilon_{1} & \mbox{if $x \in \Omega_{1}$},\\
   \varepsilon_{2} & \mbox{if $x \in \Omega_{2}$}
  \end{array} \right.
 $$
 is a piecewise constant dielectric. We decompose the solution as $\varphi = \varphi_{s} + \varphi_{p}$, where $\varphi_{s}$ is the potential due to the sources:
 \begin{equation}
  \varphi_{s} \left( x \right) = \frac{1}{\varepsilon_{1}} \sum_{i = 1}^{n} q_{i} G \left( x, x_{i} \right),
  \label{eq:mol-elect:phi_s}
 \end{equation}
 with $G$ given by (\ref{eq:3d-laplace}), and $\varphi_{p}$ is a piecewise harmonic potential, satisfying the jump conditions
 $$
  \left[ \varphi_{p} \right] = 0, \qquad \left[ \varepsilon \frac{\partial \varphi_{p}}{\partial \nu} \right] = -\left[ \varepsilon \frac{\partial \varphi_{s}}{\partial \nu} \right]
 $$
 on $\Sigma$, where $\nu$ is the unit outer normal. We can write $\varphi_{p}$, called the {\em polarization response}, as a single-layer potential \cite{guenther:1988:prentice-hall}
 \begin{equation}
  \varphi_{p} \left( x \right) = \int_{\Sigma} G \left( x, y \right) \sigma \left( y \right) dy,
  \label{eq:mol-elect:phi_p}
 \end{equation}
 which yields the boundary integral equation
 $$
  \frac{1}{2} \sigma \left( x \right) + \lambda \int_{\Sigma} \frac{\partial G}{\partial \nu_{x}} \left( x, y \right) \sigma \left( y \right) dy = -\lambda \frac{\partial \varphi_{s}}{\partial \nu} \left( x \right),
 $$
 where $\lambda = (\varepsilon_{1} - \varepsilon_{2}) / (\varepsilon_{1} + \varepsilon_{2})$, in terms of the {\em polarization charge} $\sigma$. Once $\sigma$ has been computed, the potential at any point can be evaluated using (\ref{eq:mol-elect:phi_s}) and (\ref{eq:mol-elect:phi_p}).

\begin{figure}
  \begin{center}
    \includegraphics{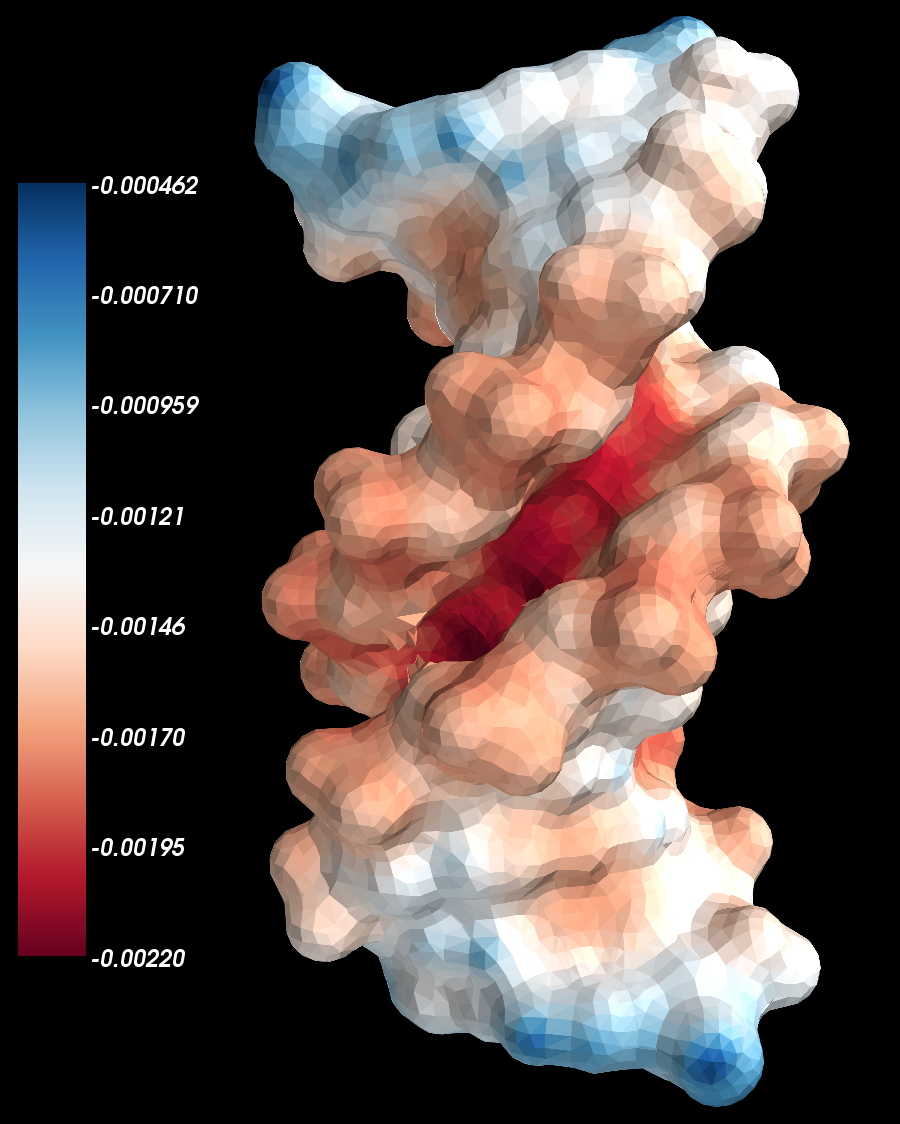}
  \end{center}
  \caption{Surface potential of DNA (PDB ID: 1BNA) in units of the elementary charge, computed using recursive skeletonization to precision $\epsilon = 10^{-3}$. The molecular surface was discretized using $N = 19752$ triangles.}
  \label{fig:mol-elect}
 \end{figure}

 We generated a molecular surface for a short segment of DNA \cite[PDB ID: 1BNA]{drew:1981:proc-natl-acad-sci-usa} consisting of $N = 19752$ triangles using MSMS \cite{sanner:1996:biopolymers}. Strengths were assigned to each of $n = 486$ heavy atoms using Amber partial charges \cite{case:2005:j-comput-chem} through PDB2PQR \cite{dolinsky:2004:nucleic-acids-res}. The system was solved with $\varepsilon_{1} = 20$ and $\varepsilon_{2} = 80$ at precision $\epsilon = 10^{-3}$; the resulting potential $\varphi$ on $\Sigma$ is shown in Fig.\ \ref{fig:mol-elect}.
 The net solution time was $T_{\RS} = 592$ s, with an inverse application time of $T_{\sv} = 0.08$ s, to be compared with $T_{\FMM} = 27$ s using FMM/GMRES. Thus, when sampling the electrostatic field for many charge configurations $\{ q_{i} \}$, as is common in computational chemistry (e.g., \cite{beroza:1991:proc-natl-acad-sci-usa}), our solver can provide a speedup provided that the number of such configurations is greater than $\sim 25$. We remark that the evaluation of $\varphi$ at fixed points, e.g., on $\Sigma$, via (\ref{eq:mol-elect:phi_s}) and (\ref{eq:mol-elect:phi_p}) can also be accelerated using our algorithm in its capacity as a generalized FMM; the computation time for this would be similar to $T_{\sv}$.

 \subsection{Multiple scattering}
 As a final example, we show how direct solvers can be combined with FMM-based iterative methods to great effect in the context of a multiple scattering problem. For this, let $\Omega_{i}$, for $i = 1, \dots, p$, be a collection of acoustic scatterers in 2D with boundaries $\Sigma_{i}$. Then, using the language of acoustics, the pressure field satisfies
 \begin{equation}
  \left( \Delta + k^{2} \right) u = 0 \quad \mbox{in $\displaystyle \mathbb{R}^{2} \setminus \bigcup_{i = 1}^{p} \Omega_{i}$}.
  \label{eq:multi-scat-helm}
 \end{equation}
 Assuming that the obstacles are {\em sound-hard}, we must compute the exterior solution that satisfies the Neumann boundary condition
 $$
  \frac{\partial u}{\partial \nu} = 0 \quad \mbox{on $\displaystyle \bigcup_{i = 1}^{p} \Sigma_{i}$}.
 $$
 If $u = u_{i} + u_{s}$, where $u_{i}$ is an incoming field satisfying (\ref{eq:multi-scat-helm}), then the scattered field $u_{s}$ also satisfies (\ref{eq:multi-scat-helm}) with boundary condition
 $$
  \frac{\partial u_{s}}{\partial \nu} = -\frac{\partial u_{i}}{\partial \nu} \quad \mbox{on $\displaystyle \bigcup_{i = 1}^{p} \Sigma_{i}$}.
 $$

 We write the scattered field as $u_{s} = \sum_{i = 1}^{p} u_{s, i}$, where
 $$
  u_{s, i} \left( x \right) = \int_{\Sigma_{i}} G \left( x, y \right) \sigma_{i} \left( y \right) dy,
 $$
 where $G$ is the single-layer kernel (\ref{eq:2d-helmholtz}). Imposing the boundary condition yields the second-kind integral equation
 $$
  -\frac{1}{2} \sigma_{i} + \sum_{j = 1}^{p} K_{ij} \sigma_{j} = -\left. \frac{\partial u_{i}}{\partial \nu} \right|_{\Sigma_{i}} \quad \mbox{on $\Sigma_{i}$}, \quad i = 1, \dots, p,
 $$
 where
 $$
  K_{ij} \sigma_{j} \left( x \right) = \int_{\Sigma_{j}} \frac{\partial G}{\partial \nu_{x}} \left( x, y \right) \sigma_{j} \left( y \right) dy \quad \mbox{for $x \in \Sigma_{i}$}.
 $$
 In operator notation, the linear system therefore has the form
 $$
  \sum_{i=1}^p A_{ij} \sigma_{j} = -\left. \frac{\partial u_{i}}{\partial \nu} \right|_{\Sigma_{i}}, \quad A_{ij} = \left\{
  \begin{array}{ll}
   -\frac{1}{2} I + K_{ii} & \mbox{if $i = j$},\\
   K_{ij} & \mbox{if $i \neq j$}.
  \end{array} \right.
 $$

 We solve this system using FMM/GMRES with the block diagonal preconditioner
 $$
  P^{-1} = \left[
  \begin{array}{ccc}
   A_{11}^{-1}\\
   & \ddots\\
   & & A_{pp}^{-1}
  \end{array} \right],
 $$
 where each $A_{ii}^{-1}$ is computed using recursive skeletonization; observe that $A_{ii}^{-1}$ is precisely the solution operator for scatterer $\Omega_{i}$ in isolation. The question is whether this preconditioner will significantly reduce the iteration count required, which is typically quite high for problems of appreciable size.

 As a test, we embedded two identical scatterers, each described in polar coordinates by the radial function $r = [2 + \cos (3 \theta)] / 6$, where $\theta$ is the polar angle; each scatterer is smooth, though somewhat complicated, and was taken to be ten wavelengths in size. We assumed an incoming field given by the plane wave $u_{i} = e^{\imath k x_{2}}$, where $x = (x_{1}, x_{2})$, and considered the scattering problem at various horizontal separation distances $\delta$ between the centers of the scatterers. Each configuration was solved both with and without the preconditioner $P^{-1}$ to precision $\epsilon = 10^{-6}$; each scatterer was discretized using a corrected trapezoidal rule \cite{kapur:1997:siam-j-numer-anal} with $N = 1024$ points.

\begin{figure}
  \begin{center}
    \includegraphics{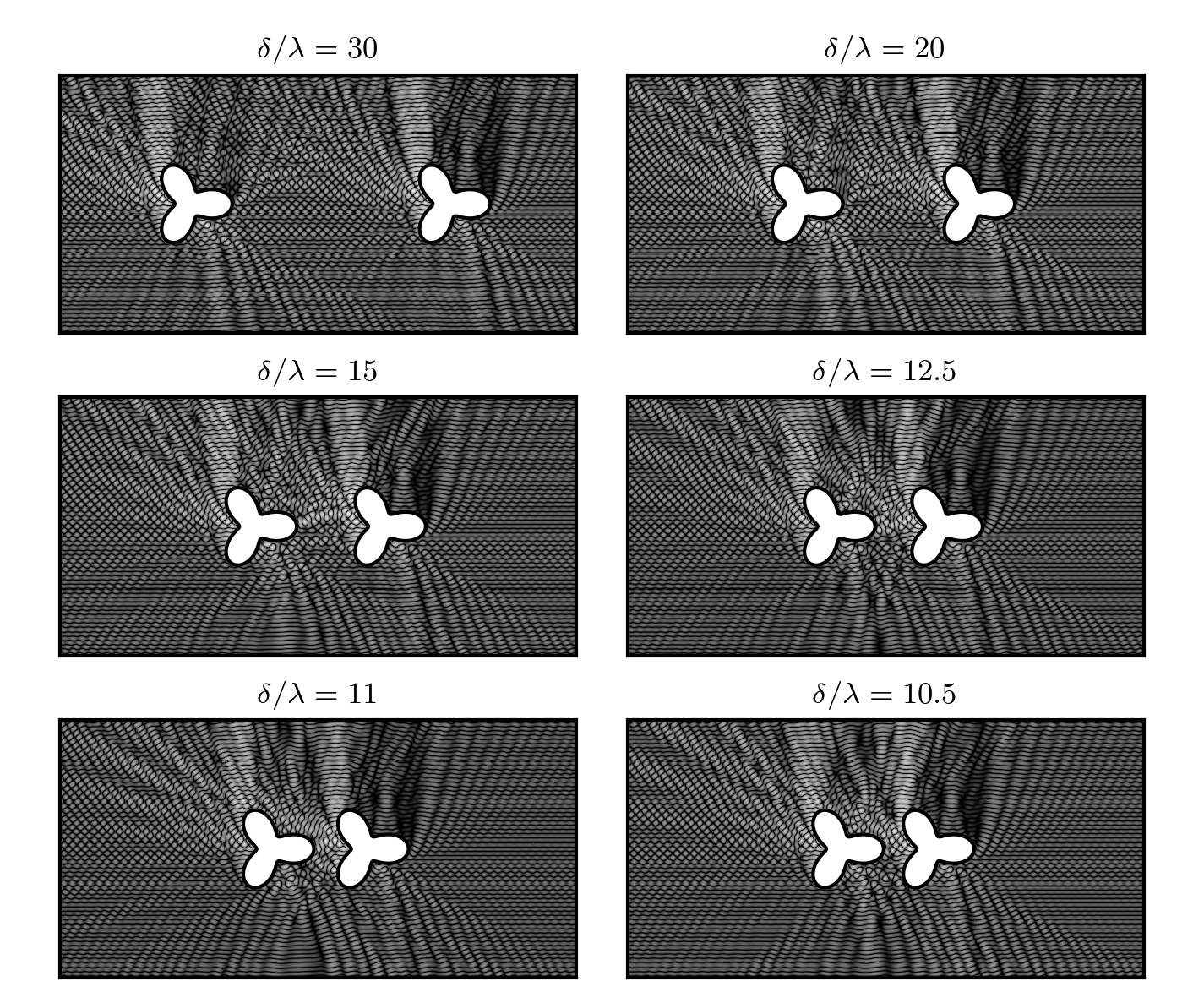}
  \end{center}
  \caption{Instantaneous intensity $[\Re (u)]^{2}$ of the pressure field in response to an incoming vertical plane wave for various scattering geometries characterized by the separation distance $\delta / \lambda$ in wavelengths between the centers of two identical scatterers.}
  \label{fig:multi-scat}
 \end{figure}

 The intensities of the resulting pressure fields are shown in Fig.\ \ref{fig:multi-scat}, with numerical data given in Table \ref{tab:multi-scat}.
 It is clear that the preconditioner is highly effective: following a precomputation time of $0.76$ s to construct $P^{-1}$, which is amortized over all solves, the number of iterations required was decreased from $n_{\FMM} \sim 700$ to just $n_{\RS} \sim 10$ for each case.
 \begin{table}
  \caption{Numerical results for the multiple scattering example, consisting of six configurations with various separation distances $\delta / \lambda$, relative to the wavelength, between the centers of two identical scatterers, solved to precision $\epsilon = 10^{-6}$: $T_{\FMM}$, time for FMM/GMRES solve (s); $T_{\RS}$, time for preconditioned FMM/GMRES solve (s); $n_{\FMM}$, number of iterations required for FMM/GMRES; $n_{\RS}$, number of iterations required for preconditioned FMM/GMRES; $E$, relative error; $T_{\cm}$, matrix compression time for scatterer (s);
$T_{\lu}$, sparse matrix factorization time for scatterer (s).}
  \label{tab:multi-scat}
  \begin{center}
   {\footnotesize
    \begin{tabular}{r|ccccc}
     \hline
     \multicolumn{1}{c|}{$\delta / \lambda$} & $T_{\FMM}$ & $T_{\RS}$ & $n_{\FMM}$ & $n_{\RS}$ & $E$\\
     \hline
     30.0 & 7.9\E$+$1 & 8.9\E$-$1 & 697 & \phantom{0}8 & 1.3\E$-$8\\
     20.0 & 7.7\E$+$1 & 1.1\E$+$0 & 694 & 10 & 5.8\E$-$9\\
     15.0 & 8.0\E$+$1 & 1.2\E$+$0 & 695 & 11 & 6.9\E$-$9\\
     12.5 & 7.9\E$+$1 & 1.3\E$+$0 & 695 & 12 & 7.8\E$-$9\\
     11.0 & 7.9\E$+$1 & 1.4\E$+$0 & 704 & 14 & 8.7\E$-$9\\
     10.5 & 8.0\E$+$1 & 1.5\E$+$0 & 706 & 14 & 1.3\E$-$8\\
     \hline
     \multicolumn{1}{c|}{$T_{\cm}$} & & 6.6\E$-$1\\
     \multicolumn{1}{c|}{$T_{\lu}$} & & 9.3\E$-$2\\
     \hline
     \multicolumn{1}{c|}{total} & 4.7\E$+$2 & 8.1\E$+$0\\
     \hline
    \end{tabular}
   }
  \end{center}
 \end{table}
 As expected, more iterations were necessary for smaller $\delta$, though the difference was not too dramatic. The ratio of the total solution time required for all solves was $\sim 60$ for the unpreconditioned versus the preconditioned method.

 \section{Generalizations and conclusions}
 \label{sec:generalizations-conclusions}
 We have presented a multilevel matrix compression algorithm and demonstrated its efficiency at accelerating matrix-vector multiplication and matrix inversion in a variety of contexts. The matrix structure required is fairly general and relies only on the assumption that the matrix have low-rank off-diagonal blocks. As a fast direct solver for the boundary integral equations of potential theory, we found our algorithm to be competitive with fast iterative methods based on FMM/GMRES in both 2D and 3D, provided that the integral equation kernel is not too oscillatory, and that the system size is not too large in 3D. In such cases, the total solution times for both methods were very comparable. Our solver has clear advantages, however, for problems with ill-conditioned matrices (in which case the number of iterations required by FMM/GMRES can increase dramatically), or those involving multiple right-hand sides (in which case the cost of matrix compression and factorization can be amortized). The latter category includes the use of our solver as a preconditioner for iterative methods, which we expect to be quite promising, particularly for large-scale 3D problems with complex geometries.

 A principal limitation of the approach described here is the growth in the cost of factorization in 3D or higher, which prohibits the scheme from achieving optimal $\bigO (N)$ or nearly optimal $\bigO(N \log N)$ complexity. It is, however, straightforward to implement and quite effective. All of the hierarchical compression-based approaches (HSS matrices \cite{chandrasekaran:2006b:siam-j-matrix-anal-appl,chandrasekaran:2006a:siam-j-matrix-anal-appl,xia:2009:siam-j-matrix-anal-appl}, $\mathcal{H}$-matrices \cite{hackbusch:1999:computing,hackbusch:2002:computing,hackbusch:2000:computing} and skeletonization \cite{gillman:2012:front-math-china,greengard:2009:acta-numer,martinsson:2005:j-comput-phys}) are capable of overcoming this obstacle. The development of simple and effective schemes that curtail this growth is an active area of research, and we expect that $\bigO (N \log N)$ direct solvers with small pre-factors in higher dimensions will be constructed shortly, at least for non-oscillatory problems. It is clear that all of these techniques provide improved solution times for high-frequency {\em volume} integral equations, due to the compression afforded by Green's theorem in moving data from the volume to the boundary. More precisely, the cost of solving high-frequency volume wave scattering problems in 2D are $\bigO (N^{3/2})$ and $\bigO (N \log N)$ for precomputation and solution, respectively. For related work, see \cite{chen:2002:adv-comput-math,winebrand:2009:intl-conf-electromag-adv-application}.

 Finally, although all numerical results have presently been reported for a single processor, the algorithm is naturally parallelizable: many computations are organized in a block sweep structure, where each block can be processed independently. This is clearly true of the recursive skeletonization phase using proxy surfaces (with a possible loss of $\bigO (\log N)$ in performance since there are $\bigO (\log N)$ levels in the hierarchy). As for the solver phase, arguments can be made in support of both the original ``hand-rolled'' Gaussian elimination approach and our framework that relies on sparse embedding.
We expect that, by making use of
UMFPACK and other state-of-the-art parallel sparse solvers (e.g., SuperLU \cite{li:2005:acm-trans-math-softw}, MUMPS \cite{amestoy:2001:siam-j-matrix-anal-appl}, Pardiso \cite{schenk:2004:future-gener-comput-syst}, WSMP \cite{gupta:2000:ibm}), our overall
strategy will help simplify the implementation of skeletonization-based schemes on
high-performance computing platforms as well.

 \section*{Acknowledgements}
 We would like to thank Zydrunas Gimbutas and Mark Tygert for many helpful discussions.

\end{document}